\documentclass[11pt,draftclsnofoot,onecolumn]{IEEEtran}
\usepackage{amsmath,amssymb,eucal,graphicx,subfig}
\usepackage{epsfig}
\usepackage{exscale}
\usepackage{pstricks}
\usepackage{pst-node}
\usepackage{pst-blur}
\usepackage{setspace}
\usepackage{url}
\usepackage{hyperref}
\definecolor{Pink}{rgb}{1.,0.75,0.8}
\DeclareMathAlphabet{\mathcal}{OMS}{cmsy}{m}{n} 

\long\def\comment#1{}

\newfont{\bbb}{msbm10 scaled 700}

\newfont{\bb}{msbm10 scaled 1100}
\newcommand{\CC}{\mbox{\bb C}}

\newcommand{\EE}{{\mathbb E}}

\newcommand{\bv}{{\bf b}}

\newcommand{\nv}{{\bf n}}

\newcommand{\wv}{{\bf w}}

\newcommand{\xv}{{\bf x}}
\newcommand{\yv}{{\bf y}}

\newcommand{\zerov}{{\bf 0}}

\newcommand{\Am}{{\bf A}}

\newcommand{\Fm}{{\bf F}}

\newcommand{\Id}{{\bf I}}

\newcommand{\Om}{{\bf O}}
\newcommand{\Pm}{{\bf P}}

\newcommand{\Sm}{{\bf S}}

\newcommand{\Xm}{{\bf X}}

\newcommand{\Cc}{{\cal C}}

\newcommand{\Nc}{{\cal N}}
\newcommand{\Oc}{{\cal O}}

\newcommand{\Sc}{{\cal S}}

\newcommand{\Zc}{{\cal Z}}

\newcommand{\Lambdam}{\hbox{\boldmath$\Lambda$}}

\newcommand{\Sigmam}{\hbox{\boldmath$\Sigma$}}

\newcommand{\Psim}{\hbox{\boldmath$\Psi$}}
\newcommand{\psim}{\hbox{\boldmath$\psi$}}

\newcommand{\omegam}{\hbox{\boldmath$\omega$}}
\newcommand{\etam}{\hbox{\boldmath$\eta$}}

\newcommand{\Thetam}{\hbox{\boldmath$\Theta$}}

\renewcommand{\det}{{\hbox{det}}}

\renewcommand{\arg}{{\hbox{arg}}}

\renewcommand{\Re}{{\rm Re}}

\newcommand{\herm}{{\sf H}}

\newcommand{\transp}{{\sf T}}
\newcommand{\calL}{\mbox{${\mathcal L}$}}
\newcommand{\calO}{\mbox{${\mathcal O}$}}

\DeclareFontFamily{U}{cmfi}{}
\DeclareFontShape{U}{cmfi}{m}{n}{ <-> cmfi10 }{}
\DeclareSymbolFont{CMFI}{U}{cmfi}{m}{n}

\newcommand{\define}{\stackrel{\triangle}{=}}

\def\BibTeX{{\rm B\kern-.05em{\sc i\kern-.025em b}\kern-.08em
    T\kern-.1667em\lower.7ex\hbox{E}\kern-.125emX}}

\title{Structure-Based Bayesian Sparse Reconstruction}
\author{Ahmed A. Quadeer$^{\rm a}$ and Tareq Y. Al-Naffouri$^{\rm b}$
\\$^{\rm a}$King Fahd University of Petroleum $\&$ Minerals, Dhahran, Saudi Arabia
\\$^{\rm b}$King Abdullah University of Science $\&$ Technology, Thuwal, Saudi Arabia
\\Email: aquadeer@kfupm.edu.sa and tareq.alnaffouri@kaust.edu.sa}
\date{\today}

\begin{document}
\doublespace
\maketitle

\begin{abstract}
 Sparse signal reconstruction algorithms have attracted research attention due to their wide applications in various fields. In this paper, we present a simple Bayesian approach that utilizes the sparsity constraint and {\em a priori} statistical information (Gaussian or otherwise) to obtain near optimal estimates. In addition, we make use of the rich structure of the sensing matrix encountered in many signal processing applications to develop a fast sparse recovery algorithm. The computational complexity of the proposed algorithm is relatively low compared with the widely used convex relaxation methods as well as greedy matching pursuit techniques, especially at a low sparsity rate.\footnote{This work was partially supported by SABIC through an internally funded project from DSR, KFUPM (Project No. SB101006) and partially by King Abdulaziz City for Science and Technology (KACST) through the Science \& Technology Unit at KFUPM (Project No. 09-ELE763-04) as part of the National Science, Technology and Innovation Plan. The work of Tareq Y. Al-Naffouri was also supported by the Fullbright Scholar Program. Part of this work was presented at the Allerton Conference on Communications, Control and Computing, USA.}
\end{abstract}

\section{Introduction}
\label{sec:intro}
Compressive Sensing/Compressed Sampling (CS) is a fairly new field of research that is finding many applications in statistics and signal processing \cite{candes-mag}. As its name suggests, CS attempts to acquire a signal (inherently sparse in some subspace) at a compressed rate by randomly projecting it onto a subspace that is much smaller than the dimension of the signal itself. Provided that the sensing matrix satisfies a few conditions, the sparsity pattern of such a signal can be recovered non-combinatorially with high probability. This is in direct contrast to the traditional approach of sampling signals according to the Nyquist theorem and then discarding the insignificant samples.
Generally, most naturally occurring signals are sparse in some basis/domain and CS can therefore be utilized for their reconstruction. CS has been used successfully in, for example (but not limited to), peak-to-average power ratio reduction in orthogonal frequency division multiplexing (OFDM) \cite{PAPR}, image processing (one-pixel camera \cite{image}), impulse noise estimation and cancellation in power-line communication and digital subscriber lines (DSL) \cite{impulse-noise}, magnetic resonance imaging (MRI) \cite{MRI}, channel estimation in communications systems \cite{ce}, ultra-wideband (UWB) channel estimation \cite{ce1}, direction-of-arrival (DOA) estimation \cite{geert}, and radar design \cite{radar}, to name a few.

The CS problem can be set up as follows. Let $\xv\in{\CC^N}$ be a $P$-sparse signal (i.e., a signal that consists of $P$ non-zero coefficients in an $N$-dimensional space with $P << N$) in some domain and let $\yv\in{\CC^M}$ be the observation vector with $M << N$ given by
\begin{equation} \label{model}
\yv = \Psim \xv + \nv
\end{equation}
where $\Psim$ is an $M \times N$ measurement/sensing matrix that is assumed to be incoherent with the domain in which $\xv$ is sparse and $\nv$ is complex additive white Gaussian noise, $\Cc\Nc(\zerov,\sigma_n^2\Id_M)$. As $M << N$, this is an ill-posed problem as there is an infinite number of solutions for  $\xv$ satisfying~(\ref{model}). Now if it is known {\em a priori} that $\xv$ is sparse, the theoretical way to reconstruct the signal is to solve an $\ell_0$-norm minimization problem using only $M = 2P$ measurements when the signal and measurements are free of noise \cite{Baraniuk}
\begin{equation}
\hat{\xv}=\min_{\xv} \| \xv \|_0 \quad \quad \mbox{subject to} \; \; \yv = \Psim \xv.
\end{equation}
Unfortunately, solving the $\ell_0$-norm minimization problem is NP-hard \cite{Baraniuk} \cite{candes-romberg-tao} and is therefore not practical. Thus, different sub-optimal approaches, categorized as compressive sensing, have been presented in the literature to solve this problem. 
In \cite{candes-romberg-tao} and \cite{candes-randall}, it has been shown that $\xv$ can be reconstructed with high probability in polynomial time by using convex relaxation approaches at the cost of an increase in the required number of measurements. 
This is done by solving a relaxed $\ell_1$-norm minimization problem using linear programming instead of $\ell_0$-norm minimization \cite{candes-romberg-tao}, \cite{candes-randall}
\begin{eqnarray} \label{a1}
\hat{\xv}=  \min_{\xv} \| \xv \|_1 \quad \quad \mbox{subject to} \; \; \| \yv - \Psim \xv \|_2 \leq \epsilon
\end{eqnarray}
where $\epsilon = \sqrt{\sigma_n^2(M+\sqrt{2M})}$. For $\ell_1$-norm minimization to reconstruct the sparse signal accurately,
the sensing matrix, $\Psim$, should be sufficiently incoherent. In other words, the coherence, defined as $\mu(\Psim) \define \max_{i \neq j} |\langle \psim_i \psim_j \rangle|$, should be as small as possible (with $\mu(\Psim) = 1$ depicting the worst case) \cite{candes-romberg-tao}.
In \cite{BCS}, it has been shown that these convex relaxation approaches have a Bayesian rendition and may be viewed as maximizing the maximum {\em a posteriori} estimate of $\xv$, given that $\xv$ has a Laplacian distribution. Although convex relaxation approaches are able to recover sparse signals by solving  under-determined systems of equations, they also suffer from a number of drawbacks (some of which are common to other sparse recovery algorithms including \cite{OMP}-\cite{COSAMP}) that we discuss below.
\subsection{Drawbacks of Convex Relaxation Approaches}

\subsubsection{Complexity}
Convex relaxation relies on linear programming to solve the convex $\ell_1$-norm minimization problem, which is computationally relatively complex (its complexity is of the order $\Oc(M^2N^{3/2})$ when interior point methods are used \cite{interior_point}). This approach can therefore not be used in problems with very large dimensions. To overcome this drawback, many greedy algorithms have been proposed that recover the sparse signal iteratively. These include Orthogonal Matching Pursuit (OMP) \cite{OMP2}, \cite{OMP}, Regularized Orthogonal Matching Pursuit (ROMP) \cite{ROMP}, Stagewise Orthogonal Matching Pursuit (StOMP) \cite{STOMP}, and Compressive Sampling Matching Pursuit (CoSamp) \cite{COSAMP}. These greedy approaches are relatively faster than their convex relaxation counterparts (approximately $\Oc(MNR)$ where $R$ is the number of iterations).

\subsubsection{The need for randomness in the sensing matrix}
Convex relaxation methods cannot make use of the structure exhibited by the sensing matrix (e.g., a structure that comes from a Toeplitz sensing matrix or that of a partial discrete Fourier transform (DFT) matrix). In fact, if anything, this structure is harmful to these methods as the best results are obtained when the sensing matrix is close to random. This comes in contrast to current digital signal processing architectures that only deal with uniform sampling. We would thus like to employ more feasible and standard sub-sampling approaches. 

\subsubsection{Inability to harness {\em a priori} statistical information}
Convex relaxation methods are not able to take account of any {\em a priori} statistical information (apart from sparsity information) about the signal support and additive noise.
Any {\em a priori} statistical information can be used on the result obtained from the convex relaxation method to refine both the signal support obtained and the resulting estimate through a hypothesis testing approach \cite{ISCAS}. However, this is only useful if these approaches are indeed able to recover the signal support. In other words, performance is bottle-necked by the support recovering capability of these approaches. We note here that the use of {\em a priori} statistical information for sparse signal recovery has been studied in a Bayesian context in \cite{BCS} and in algorithms based on belief propagation \cite{Tan}, \cite{Montanari}.
Both \cite{FBMP} and \cite{Larsson} use {\em a priori} statistical information (assuming $\xv$ to be mixed Bernoulli-Gaussian); only \cite{FBMP} uses this information in a recursive manner to obtain a fast sparse signal recovery algorithm. However, it is not clear how these approaches can be extended to the non-Gaussian case.

\subsubsection{Evaluating performance in statistically familiar terms}
It is difficult to quantify the performance of convex relaxation estimates analytically in terms of the mean squared error (MSE) or bias or to relate these estimates to those obtained through more conventional approaches, e.g., maximum a posteriori probability (MAP), minimum mean-square error (MMSE), or maximum likelihood (ML).\footnote{It is worth noting that convex relaxation approaches have their merit in that they are agnostic to the signal distribution and thus can be quite useful when worst-case analysis is desired as opposed to average-case analysis.}

\subsubsection{Trading performance for computational complexity}
In general, convex relaxation approaches do not exhibit the customary tradeoff between increased computational complexity and improved recovery as is the case for, say, iterative decoding or joint channel and data detection. Rather, they solve some $\ell_1$ problem using (second-order cone programming) with a set complexity. A number of works have attempted to derive sharp thresholds for support recovery \cite{wainwright}, \cite{goyal}. In other words, the only degree of freedom available for the designer to improve performance is to increase the number of measurements. Several iterative implementations \cite{yin}, \cite{SPGL1} of convex relaxation approaches provide some sort of flexibility by trading performance for complexity.


\subsection{Motivation and Paper Organization}

In this paper, we present a Bayesian approach to sparse signal recovery that has low complexity and makes a collective use of 1) {\em a priori} statistical properties of the signal and noise, 2) sparsity information, and 3) the rich structure of the sensing matrix, $\Psim$. Although there have been some works that use the structure of the sensing matrix (e.g., \cite{tree}), it has not yet been rigorously exploited to aid in algorithm development and complexity reduction.
We also show how our approach is able to deal with both Gaussian and non-Gaussian (or unknown) priors, and how we can compute performance measures of our estimates. In essence, we demonstrate how our technique enables us to tackle all the drawbacks of convex relaxation approaches.

This remainder of this paper is organized as follows. We start by describing the signal model in the next section. In Section \ref{sec:optimum-estimate}, we derive the MMSE/MAP estimates and introduce the various terms that need to be evaluated. In Section \ref{sec:cont-band}, we demonstrate how the structure of the sensing matrix can be used to recover the sparse signal in a divide-and-conquer manner. Section \ref{orth_algo} details the proposed sparse reconstruction algorithm that we call Orthogonal Clustering. 
Section~\ref{sec:complexity} presents the different structural properties of the sensing matrix that are exploited by the proposed algorithm to reduce the computational complexity.
The performance of the proposed algorithm is compared with various sparse reconstruction algorithms presented in the literature by numerical simulations in Section \ref{sec:sim}, which is followed by our conclusions in Section~\ref{sec:conc}.

\subsection{Notation}

We denote scalars with lower-case letters (e.g., $x$), vectors with lower-case bold-faced letters (e.g., $\xv$), matrices with upper-case, bold-faced letters (e.g., $\Xm$), and sets with script notation (e.g. $\Sc$). We use $\xv_i$ to denote the $i^{th}$ column of matrix $\Xm$, $x(j)$ to denote the $j^{th}$ entry of vector $\xv$, and $\Sc_i$ to denote a subset of a set $\Sc$. We also use $\Xm_{\Sc}$ to denote the sub-matrix formed by the columns $\{\xv_i : i \in \Sc \}$, indexed by the set $\Sc$. Finally, we use $\hat{\xv}$, $\xv^*$, $\xv^\transp$, and $\xv^\herm$ to respectively denote the estimate, conjugate, transpose, and conjugate transpose of a vector $\xv$.


\section{Signal Model}\label{sec:problem-formulation}


We adopt the signal model in (\ref{model}). Here, the vector $\xv$ is modelled as  
$\xv = \xv_B \odot \xv_G$,
where $\odot$ denotes the Hadamard (element-by-element) multiplication. The entries of $\xv_B$ are independent and identically distributed (i.i.d) Bernoulli random variables and the entries of $\xv_G$ are drawn identically and independently from some zero mean distribution.\footnote{Most of the results presented in this paper also apply to the case when the entries are independent but not necessarily identically distributed.} In other words, we assume that $x_B(i)$s are Bernoulli with success probability $p$ and similarly that the $x_G(i)$s are i.i.d variables with marginal probability distribution function $f(x)$. The noise $\nv$ is assumed to be complex circularly symmetric Gaussian, i.e., $\nv \sim \Cc\Nc(0,\sigma_n^2 \Id_M)$.
When the support set $\Sc$ of $\xv$ is known, we can equivalently write (\ref{model}) as
\begin{equation} \label{model-S}
\yv = \Psim_\Sc \xv_\Sc + \nv.
\end{equation}

\section{Optimum Estimation of $\xv$} \label{sec:optimum-estimate}
Our task is to obtain the optimum estimate of $\xv$ given the observation $\yv$. We can pursue either an MMSE or a MAP approach to achieve this goal. In the following, we elaborate on how we can obtain these two estimates.

\subsection{MMSE Estimation of $\xv$}
The MMSE estimate of $\xv$ given the observation $\yv$ can be expressed as
\begin{equation} \label{x_mmse}
\hat{\xv}_{\rm MMSE} = \EE[\xv|\yv] = \sum_{\Sc} p(\Sc|\yv) \EE[\xv|\yv,\Sc]
\end{equation}
where the sum is over all the possible support sets $\Sc$ of $\xv$. 
The likelihood and expectation involved in (\ref{x_mmse}) are evaluated below.

\subsubsection{Evaluation of $\EE[\xv|\yv,\Sc]$}
Recall that the relationship between $\yv$ and $\xv$ is linear (see (\ref{model})). Thus, in the case when $\xv$ conditioned on
its support is Gaussian, $\EE[\xv|\yv,\Sc]$ is nothing but the linear MMSE estimate of $\xv$
given $\yv$ (and $\Sc$), i.e.,
\begin{equation}\label{Exp_MMSE}
\EE[\xv_\Sc|\yv] \define \EE[\xv|\yv,\Sc] = \sigma_x^2 \Psim^\herm_{\Sc} \Sigmam_\Sc^{-1}\yv
\end{equation}
where 
\begin{equation}\label{cov-yprime}
\Sigmam_\Sc = \frac{1}{\sigma_n^2} \EE[ \yv \yv^\herm |\Sc ] =  \Id_M +
\frac{\sigma_x^2}{\sigma_n^2} \Psim_{\Sc} \Psim_{\Sc}^\herm.
\end{equation}
When $\xv|{\Sc}$ is non-Gaussian or when its statistics are unknown, the expectation $\EE[\xv|\yv,\Sc]$ is difficult or even impossible to calculate. Thus, we replace it by the best linear unbiased estimate (BLUE), i.e.,
\begin{equation}\label{Exp_LS}
\EE[\xv_\Sc|\yv] = (\Psim^\herm_{\Sc} \Psim_{\Sc})^{-1}\Psim^\herm_{\Sc}\yv.
\end{equation}

\subsubsection{Evaluation of $p(\Sc|\yv)$} \label{sec:Ie_prob}
Using Bayes' rule, we can rewrite $p(\Sc|\yv)$ as
\begin{equation}\label{Ie_prob}
p(\Sc|\yv) = \frac{p(\yv|\Sc)p(\Sc)}{\sum_{\Sc}p(\yv|\Sc)p(\Sc)}.
\end{equation}
As the denominator $\sum_{\Sc}p(\yv|\Sc)p(\Sc)$ is common to all posterior likelihoods, $p(\Sc|\yv)$, it is a normalizing constant that can be ignored. To evaluate $p(\Sc)$, note that the elements of $\xv$ are active according to a Bernoulli process with
success probability $p$. Thus, $p(\Sc)$ is given by
\begin{equation}\label{pSc}
p(\Sc) = p^{|\Sc|} ( 1 - p)^{N-|\Sc|}.
\end{equation}
It remains to evaluate $p(\yv|\Sc)$. Here, we distinguish between the cases of whether or not $\xv|\Sc$ is Gaussian.
\paragraph*{\em 1. $\xv|\Sc$ is Gaussian}
When $\xv|\Sc$ is Gaussian, $\yv$ is Gaussian too with zero mean and covariance $\Sigmam_\Sc$
and we can write the likelihood function as\footnote{$\|\bv \|^2_{\Am} \stackrel{\Delta}{=} \bv^\herm\Am\bv$}
\begin{equation} \label{map-metric}
p(\yv | \Sc) = \frac{\exp
\left ( - \frac{1}{\sigma_n^2} \| \yv \|^2_{\Sigmam_\Sc^{-1}} \right )
}{ \det \left ( \Sigmam_\Sc \right ) }
\end{equation}
up to an irrelevant constant multiplicative factor, ($\frac{1}{\pi^M}$).

\paragraph*{\em 2. $\xv|\Sc$ is non-Gaussian or unknown}
Alternatively, we can treat $\xv$ as a random vector of unknown (non-Gaussian) distribution, with
support $\Sc$. Therefore, given the support $\Sc$, all we
can say about $\yv$ is that it is formed by a vector in the
subspace spanned by the columns of $\Psim_\Sc$,  plus a white
Gaussian noise vector, $\nv$. It is difficult to quantify the distribution of $\yv$ even if we know the distribution of (the non-Gaussian) $\xv$. One way around this is to annihilate the non-Gaussian component and retain the Gaussian one. We do so by projecting $\yv$ onto the orthogonal complement of the span of the columns of $\Psim_\Sc$, i.e., multiplying $\yv$ by
$\Pm_\Sc^\perp = \Id - \Psim_{\Sc} \left ( \Psim_{\Sc}^\herm \Psim_{\Sc} \right )^{-1} \Psim_{\Sc}^\herm.$
This leaves us with $\Pm_\Sc^\perp \yv = \Pm_\Sc^\perp \nv$, which is zero mean and with covariance $\Pm_\Sc^\perp \sigma_n^2 \Pm_\Sc^{\perp^\herm} = \sigma_n^2 \Pm_\Sc^\perp.$
Thus, the conditional density of $\yv$ given $\Sc$ is approximately given by 
\begin{eqnarray}
\label{map-metric1} p(\yv | \Sc) \simeq \exp \left( -\frac{1}{\sigma_n^2} \left\| \Pm_\Sc^\perp \yv \right\|^2 \right).
\end{eqnarray}


\subsection{MAP Estimation of $\xv$}
To obtain the MAP estimate of $\xv$, we first determine the MAP estimate of $\Sc$, which is given by
\begin{equation}\label{S_map}
\hat{\Sc}_{\rm MAP} = \arg \max_\Sc p(\yv|\Sc) p(\Sc).
\end{equation}
The prior likelihood $p(\yv|\Sc)$, is given by (\ref{map-metric}) when $\xv|\Sc$ is Gaussian and by (\ref{map-metric1}) when $\xv|\Sc$ is non-Gaussian or unknown, whereas $p(\Sc)$ is evaluated using (\ref{pSc}). The maximization is performed over all possible $2^N$ support sets. The corresponding MAP estimate of $\xv$ is given by
\begin{equation}
\hat{\xv}_{\rm MAP} = \EE[\xv|\yv,\hat{\Sc}_{\rm MAP}].
\end{equation}
One can easily see that the MAP estimate is a special case of the MMSE estimate in which the sum (\ref{x_mmse}) is reduced to one term. As a result, we confine the discussion in the rest of the paper to MMSE estimation.

\subsection{Evaluation over $\Sc$}\label{eval-over-Sc}
Having evaluated the posterior probability and expectation, it remains to evaluate this over $2^N$ possible supports (see (\ref{x_mmse}) and (\ref{S_map})) which is a computationally daunting task. This is compounded by the fact that the calculations required for each support set in $\Sc$ are relatively expensive, requiring some form of matrix multiplication/inversion as can be seen from (\ref{Exp_MMSE})-(\ref{map-metric1}).
One way around this exhaustive approach is somehow to guess at a superset $\Sc_r$ consisting of the most probable support and limit the sum in (\ref{x_mmse}) to the superset $\Sc_r$ and its subsets, reducing the evaluation space to $2^{|\Sc_r|}$ points.
There are two techniques that help us guess at such a set $\Sc_r$.
\paragraph*{\em 1. Convex Relaxation}
Starting from (\ref{model}), we can use the standard convex relaxation tools  \cite{candes-romberg-tao}, \cite{candes-randall} to find the most probable support set, $\Sc_r$, of the sparse vector $\xv$. This is done by solving (\ref{a1}) and retaining some largest $P$ non-zero values where $P$ is selected such that P$( \left\| \Sc \right\|_0 > P)$ is very small.\footnote{As $\left\| \Sc \right\|_0$ is a binomial distribution $\sim$ B$(N,p)$, it can be approximated by a Gaussian distribution $\sim$ $\Nc(Np,Np(1-p))$, when $Np > 5$ (the DeMoivre-Laplace approximation \cite{Stoch}). In this case, P$( \left\| \Sc \right\|_0 > P) = \frac{1}{2} {\rm{erfc}} \left( \frac{P-N(1-p)}{\sqrt{(2Np(1-p))}} \right)$.}
\paragraph*{\em 2. Fast Bayesian Matching Pursuit (FBMP)}
A fast Bayesian recursive algorithm is presented in \cite{FBMP} that determines the dominant support and
the corresponding MMSE estimate of the sparse vector.\footnote{FBMP applies to the Bernoulli Gaussian case only.}
It uses a greedy tree search over all combinations in pursuit of the dominant supports. The algorithm
starts with zero active element support set. At each step, an active element is added that maximizes the
Gaussian log-likelihood function similar to (\ref{map-metric}). This procedure is repeated until we reach $P$ active elements in a branch. The procedure creates $D$ such branches, which represent a tradeoff between performance and complexity.\footnote{Though other greedy algorithms \cite{OMP}-\cite{COSAMP} can also be used, we focus here on FBMP as it utilizes {\em a priori} statistical information along with sparsity information.}

The discussion in this section applies irrespective of the type of the sensing matrix, $\Psim$. However, in many applications in signal processing and communications, the sensing matrix is highly structured. This fact, which has been largely overlooked in the CS literature, is utilized in the following to evaluate the MMSE (MAP) estimate at a much lower complexity than is currently available.
\section{A Structure-Based Bayesian Recovery Approach} \label{sec:cont-band}

Whereas in most CS literature, the sensing matrix, $\Psim$, is assumed to be drawn from a random constellation \cite{candes-romberg-tao}, \cite{candes-randall}, in many signal processing and communications applications, this matrix is highly structured. Thus, $\Psim$ could be a partial DFT matrix \cite{impulse-noise} or a Toeplitz matrix (encountered in many convolution applications \cite{ce}). Table \ref{t1} lists various possibilities of structured $\Psim$.

\begin{table}[h]
\begin{center}
\caption{Applications involving structured sensing matrices}
\begin{tabular}{|l|l|}
\hline {\scriptsize \bf Matrix $\Psim$} & {\scriptsize \bf Application} \\
\hline
\hline \scriptsize Partial DFT & \scriptsize OFDM applications including peak-to-average power ratio\\
& \scriptsize reduction \cite{PAPR}, narrow-band interference cancelation \cite{nbi},\\
& \scriptsize and impulsive noise estimation and mitigation in DSL \cite{impulse-noise}\\
\hline \scriptsize Toeplitz & \scriptsize Channel estimation \cite{ce}, UWB \cite{ce1}, and DOA estimation \cite{geert} \\
\hline \scriptsize Hankel & \scriptsize Wide-band spectrum sensing \cite{geert2}\\
\hline \scriptsize DCT & \scriptsize Image compression \cite{DCT} \\
\hline \scriptsize Structured Binary & \scriptsize Multi-user detection and contention resolution \cite{giannakis}, \cite{mu_tabish} and \\
& \scriptsize feedback reduction \cite{feedback}, \cite{hanly}\\
\hline
\end{tabular}\label{t1}
\end{center}
\end{table}

Since $\Psim$ is a fat matrix $(M << N)$, its columns are not orthogonal (in fact not even linearly independent). However, in the aforementioned applications, one can usually find an orthogonal subset of the columns of $\Psim$ that span the column space of $\Psim$. We can collect these columns into a square matrix, $\Psim_M$. The remaining $N-M$ columns of $\Psim$ group around these orthogonal columns to form semi-orthogonal clusters.
In general, the columns of $\Psim$ can be rearranged such that the farther two columns are from each other, the lower their correlation is. In this section, we demonstrate how semi-orthogonality helps to evaluate the MMSE estimate in a divide-and-conquer manner. Before we do so, we present two sensing matrices that exhibit semi-orthogonality.
\subsection{Examples of Sensing Matrices with Semi-Orthogonality}
\subsubsection{DFT Matrices}
We focus here on the case when the sensing matrix is a partial DFT matrix, i.e.,
$\Psim~=~\Sm \Fm_N$,
where $\Fm_N$ denotes the $N \times N$ unitary DFT matrix, $[\Fm_N]_{a,b} = \frac{1}{\sqrt{N}}e^{-j2\pi ab/N}$ with $a,b \in \{0, 1, \ldots,N-1\}$ and $\Sm$ is an $M \times N$ selection matrix consisting of zeros with exactly one entry equal to $1$ per row. To enforce the desired semi-orthogonal structure, the matrix $\Sm$ usually takes the form 
$\Sm = \left[ \Om_{M \times Z} \;\; \Id_{M \times M} \;\; \Om_{M \times (N-Z-M)} \right]$,
for some integer $Z$.
In other words, the sensing matrix consists of a continuous band of sensing frequencies. This is not unusual since in many OFDM problems, the band of interest (or the one free of transmission) is continuous.
In this case, the correlation between two columns can be shown to be 
\begin{eqnarray}
\label{corr_orth} \psim_k^\herm \psim_{k'} = \left\{ \begin{array}{ccc} 1, & & (k = k') \\ \left| \frac{\sin \left( \pi(k-k')M/N \right)}{M \sin \left( \pi(k-k')/N \right)} \right|,
 & & (k \neq k') \end{array} \right.
\end{eqnarray}
which is a function of the difference, $(k-k'){\hspace{-0.08in}}\mod N$. 
It thus suffices to consider the correlation of one column with the remaining ones.
Figure \ref{corr_N_DFT} illustrates this correlation for $N = 1024$ and $M = 256$.
It is worth noting that the matrix $\Psim$ exhibits other structural properties (e.g., the fact that it is a Vandermonde matrix), which helps us reduce the complexity of the MMSE estimation (see Section~\ref{sec:complexity} for further details).

\subsubsection{Toeplitz/Hankel Matrices}
We focus here on the Toeplitz case. The discussion can be easily extended to the Hankel case. A sub-sampled convolutional linear system can be written in the following matrix form,
$\yv = \Psim \xv + \nv$,
where $\yv$ is a vector of length $M$, $\xv$ is a vector of length $N$ and $\Psim$ is the $M \times N$ block Toeplitz/diagonal matrix
\[
\Psim =
\left[ \begin{array}{cccc}
\Thetam & \Om & \cdots & \Om \\
\Om & \Thetam & \cdots & \Om \\
\vdots & \ddots & \ddots & \vdots \\
\Om & \Om & \cdots & \Thetam
\end{array} \right]
\]
where the size of $\Thetam$ depends on the sub-sampling ratio. Here,
$\psim_k^\herm \psim_{k'} = 0$ for $|k - k'| > L$,
and thus the columns of $\Psim$ can easily be grouped into truly orthogonal clusters. Note also that the individual columns of $\Thetam$ are related to each other by a shift property, which we explore for further reduction in complexity in Section \ref{sec:complexity}.
\subsection{Using Orthogonality for MMSE Estimation}
Let $\Sc$ be a possible support of $\xv$. The columns of $\Psim_\Sc$ in (\ref{model-S}) can be grouped into a maximum of $C$ semi-orthogonal clusters, i.e.,
$\Psim_\Sc = [\Psim_{\Sc_1} \; \Psim_{\Sc_2} \;\cdots\; \Psim_{\Sc_C}]$, where $\Sc_i$ is the support set corresponding to the $i^{th}$ cluster (with $i = 1, 2, \cdots C$).\footnote{Here, we denote the maximum number of clusters formed by $C$ to distinguish it from $P$, that refers to the estimate of the number of active supports as in  \cite{FBMP} 
(see footnote 5). In our approach, $C$ is random and depends on a threshold. This threshold is obtained using the {\em a priori} statistical information of the noise signal, $\nv$. The procedure of forming semi-orthogonal clusters is presented in Section $\ref{orth_algo}$.} Based on this fact, (\ref{model-S}) can be written as
\begin{equation}\label{y_mmse}
\yv = \left[\Psim_{\Sc_1} \; \Psim_{\Sc_2} \;\cdots\; \Psim_{\Sc_C} \right] \left[ \begin{array}{c}  \xv_1 \\ \xv_2 \\ \vdots \\ \xv_C  \end{array} \right]
+ \nv.
\end{equation}
Columns indexed by these sets should be semi-orthogonal, i.e., $\Psim_{\Sc_i}^\herm \Psim_{\Sc_j} \simeq 0$; otherwise, $\Sc_i$ and $\Sc_j$ are merged into a bigger superset. Now, the MMSE estimate of $\xv$ simplifies to\footnote{In writing an expression like the one in (\ref{x_ammse}), it is understood that estimates of elements of $\xv$ that do not belong to $\bigcup\Sc_i$ are identically zero.}
\begin{equation}\label{x_ammse}
\hat{\xv}_{\rm MMSE} = \sum_{\Zc \subset \bigcup\Sc_i} p(\Zc|\yv) \EE[\xv|\yv,\Zc].
\end{equation}
In the following, we show that
$\hat{\xv}_{\rm MMSE}$ can be evaluated in a divide-and-conquer manner by treating each cluster independently. To do so, we present in the following how orthogonality manifests itself in the calculation of the expectation and likelihood.

\subsubsection{The effect of orthogonality on the likelihood calculation}
Recall that up to a constant factor, the likelihood can be written as $p(\Zc|\yv) = p(\yv|\Zc)p(\Zc)$. Now,
\begin{eqnarray}
\nonumber p(\Zc) &=& p(\bigcup \Zc_i) \\
\nonumber &=& p^{|\bigcup \Zc_i|} (1-p)^{N-|\bigcup \Zc_i|} \\
\nonumber &=& p^{|\Zc_1| + |\Zc_2| + \; \cdots \; + |\Zc_C|} (1-p)^{N-(|\Zc_1| + |\Zc_2| + \; \cdots \; + |\Zc_C|)}\\
 \label{pZc} &=& p(\Zc_1) p(\Zc_2) \cdots p(\Zc_C)
\end{eqnarray}
where the equality in (\ref{pZc}) is true up to some constant factor. Now, to evaluate $p(\yv | \Zc)$, we distinguish between the Gaussian and non-Gaussian cases. For brevity, we focus here on the Gaussian case and extrapolate the results to the non-Gaussian case. Recall that
\begin{eqnarray}
p(\yv | \Zc) = \frac{\exp \left ( -\frac{1}{\sigma_n^2} \| \yv \|^2_{\Sigmam_\Zc^{-1}} \right )}{ \det \left ( \Sigmam_\Zc \right ) }
\end{eqnarray}
with $\Sigmam_\Zc = \Id_M + \frac{\sigma_x^2}{\sigma_n^2} \Psim_{\Zc} \Psim_{\Zc}^\herm$. Here,
$\Psim_\Zc = \left[ \Psim_{\Zc_1} \;\; \Psim_{\Zc'} \right],$
where
$\Psim_{\Zc'} = \left[ \Psim_{\Zc_2} \;\; \Psim_{\Zc_3} \;\; \cdots \;\; \Psim_{\Zc_C} \right].$
Using the matrix inversion lemma, we can write $\Sigmam_\Zc^{-1}$ as
\begin{eqnarray}
\nonumber \Sigmam_\Zc^{-1} &=& (\Id_M + \frac{\sigma_x^2}{\sigma_n^2}\Psim_\Zc \Psim_\Zc^\herm)^{-1} = (\Id_M + \frac{\sigma_x^2}{\sigma_n^2}\Psim_{\Zc_1} \Psim_{\Zc_1}^\herm + \frac{\sigma_x^2}{\sigma_n^2}\Psim_{\Zc'} \Psim_{\Zc'}^\herm)^{-1} \\
&=& \label{cluster1} \Sigmam_{\Zc_1}^{-1} - \frac{\sigma_x^2}{\sigma_n^2} \Sigmam_{\Zc_1}^{-1} \Psim_{\Zc'} (\Id_{\Zc'} + \frac{\sigma_x^2}{\sigma_n^2} \Psim_{\Zc'}^\herm \Sigmam_{\Zc_1}^{-1} \Psim_{\Zc'})^{-1} \Psim_{\Zc'}^\herm \Sigmam_{\Zc_1}^{-1}
\end{eqnarray}
where
$\Sigmam_{\Zc_1} = \Id_M + \frac{\sigma_x^2}{\sigma_n^2} \Psim_{\Zc_1} \Psim_{\Zc_1}^\herm$.
As $\Psim_{\Zc_1}$ and $\Psim_{\Zc'}$ are almost orthogonal (i.e., $\Psim_{\Zc_1}^\herm\Psim_{\Zc'} = \Psim_{\Zc'}^\herm\Psim_{\Zc_1} \simeq 0$), (\ref{cluster1}) becomes
\begin{eqnarray}
\nonumber \Sigmam_\Zc^{-1} &=& \Id_M - \frac{\sigma_x^2}{\sigma_n^2} \Psim_{\Zc_1}( \Id_{\Zc_1} +  \frac{\sigma_x^2}{\sigma_n^2} \Psim_{\Zc_1}^\herm \Psim_{\Zc_1})^{-1}\Psim_{\Zc_1}^\herm - \frac{\sigma_x^2}{\sigma_n^2} \Psim_{\Zc'}( \Id_{\Zc'} +  \frac{\sigma_x^2}{\sigma_n^2} \Psim_{\Zc'}^\herm \Psim_{\Zc'})^{-1}\Psim_{\Zc'}^\herm\\
\nonumber &=& -\Id_M + \left(\Id_M - \frac{\sigma_x^2}{\sigma_n^2} \Psim_{\Zc_1}( \Id_{\Zc_1} +  \frac{\sigma_x^2}{\sigma_n^2} \Psim_{\Zc_1}^\herm \Psim_{\Zc_1})^{-1}\Psim_{\Zc_1}^\herm \right) + \left(\Id_M - \frac{\sigma_x^2}{\sigma_n^2} \Psim_{\Zc'}( \Id_{\Zc'} +  \frac{\sigma_x^2}{\sigma_n^2} \Psim_{\Zc'}^\herm \Psim_{\Zc'})^{-1}\Psim_{\Zc'}^\herm \right)\\
&\simeq& \label{eqn45}-\Id_M + \left(\Id_M + \frac{\sigma_x^2}{\sigma_n^2} \Psim_{\Zc_1} \Psim_{\Zc_1}^\herm \right)^{-1} + \left(\Id_M + \frac{\sigma_x^2}{\sigma_n^2} \Psim_{\Zc'} \Psim_{\Zc'}^\herm \right)^{-1}.
\end{eqnarray}
Continuing in the same manner, it is easy to show that
\begin{equation}
\Sigmam_\Zc^{-1} \simeq -(C-1)\Id_M + \sum_{i = 1}^C \left( \Id_M +  \frac{\sigma_x^2}{\sigma_n^2} \Psim_{\Zc_i} \Psim_{\Zc_i}^\herm \right)^{-1}.
\end{equation}
As such, we can write
\begin{eqnarray}
\label{cluster2} \exp\left( -\frac{1}{\sigma_n^2} \| \yv \|^2_{\Sigmam_\Zc^{-1}} \right)  \simeq \exp\left( \frac{C-1}{\sigma_n^2} \| \yv \|^2 \right) \prod_{i = 1}^C \exp\left( -\frac{1}{\sigma_n^2} \| \yv \|^2_{\Sigmam_{\Zc_i}^{-1}} \right)
\end{eqnarray}
where
$\Sigmam_{\Zc_i} = \Id_M + \frac{\sigma_x^2}{\sigma_n^2} \Psim_{\Zc_i} \Psim_{\Zc_i}^\herm.$
Using a similar procedure, we can decompose $\det(\Sigmam_\Zc)$ as
\begin{eqnarray}
\nonumber \det(\Sigmam_\Zc) &=& \det(\Id_M + \frac{\sigma_x^2}{\sigma_n^2}\Psim_{\Zc_1} \Psim_{\Zc_1}^\herm + \frac{\sigma_x^2}{\sigma_n^2}\Psim_{\Zc'} \Psim_{\Zc'}^\herm)\\
&=& \label{cluster32} \det(\Id_M + \frac{\sigma_x^2}{\sigma_n^2}\Psim_{\Zc_1} \Psim_{\Zc_1}^\herm) \det(\Id_M + \frac{\sigma_x^2}{\sigma_n^2} \Psim_{\Zc'}^\herm \Sigmam_{\Zc_1}^{-1} \Psim_{\Zc'})\\
&\simeq& \label{cluster31} \det(\Id_M + \frac{\sigma_x^2}{\sigma_n^2} \Psim_{\Zc_1} \Psim_{\Zc_1}^\herm) \det(\Id_M + \frac{\sigma_x^2}{\sigma_n^2} \Psim_{\Zc'} \Psim_{\Zc'}^\herm)\\
&=& \label{cluster33} \det(\Sigmam_{\Zc_1})\det(\Sigmam_{\Zc'})
\end{eqnarray}
where in going from (\ref{cluster32}) to (\ref{cluster31}), we used the fact that $\Psim_{\Zc_1}$ and $\Psim_{\Zc'}$ are almost orthogonal. Continuing in the same way, we can show that
\begin{eqnarray}
\label{cluster3} \det(\Sigmam_\Zc) \simeq \prod_{i = 1}^C \det(\Sigmam_{\Zc_i}). 
\end{eqnarray}
Combining (\ref{cluster2}) and (\ref{cluster3}), we obtain (up to an irrelevant multiplicative factor)
\begin{equation}\label{cluster4}
p(\yv|\Zc) \simeq \prod_{i=1}^{C} p(\yv|\Zc_i).
\end{equation}
Orthogonality allows us to reach the same conclusion (\ref{cluster4}) for the non-Gaussian case. Now, combining (\ref{pZc}) and (\ref{cluster4}), we can finally write
\begin{equation}\label{cluster5}
p(\Zc|\yv) \simeq \prod_{i=1}^C p(\Zc_i|\yv)
\end{equation}
which applies equally to the Gaussian and non-Gaussian cases.

\subsubsection{The effect of orthogonality on the expectation calculation}
In evaluating the expectation, we again distinguish between the Gaussian and non-Gaussian cases. We focus here on the non-Gaussian case for which
$\EE[\xv_\Zc|\yv] = (\Psim_\Zc^\herm \Psim_\Zc)^{-1} \Psim_\Zc^\herm \yv.$
Using the decomposition into semi-orthogonal clusters $\Psim_\Zc = [\Psim_{\Zc_1} \; \Psim_{\Zc_2} \; \cdots \Psim_{\Zc_C}]$, we can write
\begin{eqnarray}
\nonumber (\Psim_\Zc^\herm \Psim_\Zc)^{-1} \Psim_\Zc^\herm \yv &=& \left[ \begin{array}{cccc}
\Psim_{\Zc_1}^\herm \Psim_{\Zc_1} & \Psim_{\Zc_1}^\herm \Psim_{\Zc_2} & \cdots & \Psim_{\Zc_1}^\herm \Psim_{\Zc_C} \\
\vdots & \vdots & \ddots & \vdots \\
\Psim_{\Zc_C}^\herm \Psim_{\Zc_1} & \Psim_{\Zc_C}^\herm \Psim_{\Zc_2} & \cdots & \Psim_{\Zc_C}^\herm \Psim_{\Zc_C} \\
\end{array} \right]^{-1}
\left[\begin{array}{c} \Psim_{\Zc_1}^\herm \yv \\ \vdots \\ \Psim_{\Zc_C}^\herm \yv \end{array}\right] \\
\nonumber &\simeq& \left[\begin{array}{c}
 (\Psim_{\Zc_1}^\herm \Psim_{\Zc_1})^{-1} \Psim_{\Zc_1}^\herm \yv \\
 \vdots \\
 (\Psim_{\Zc_C}^\herm \Psim_{\Zc_C})^{-1} \Psim_{\Zc_C}^\herm \yv \\
\end{array} \right]\\
\mbox{i.e.,} \quad
\label{Exp_Zc}
\EE[\xv_\Zc|\yv] &\simeq&
\left[\begin{array}{c}
\EE[\xv_{\Zc_1}|\yv] \\
\vdots \\
\EE[\xv_{\Zc_C}|\yv]
\end{array} \right].
\end{eqnarray}
Orthogonality allows us to write an identical expression to (\ref{Exp_Zc}) in the Gaussian case.
\subsubsection{The effect of orthogonality on the MMSE estimation}
We are now ready to show how (semi)orthogonality helps with the MMSE evaluation. To do this, we substitute the decomposed expressions (\ref{cluster5}) and (\ref{Exp_Zc}) into (\ref{x_ammse}) to get
\begin{eqnarray}
\nonumber \hat{\xv}_{\rm MMSE} &=& \sum_{\Zc \subset \bigcup\Sc_i} p(\Zc|\yv) \EE[\xv|\yv,\Zc]\\
\nonumber &\simeq& \sum_{\Zc_i \subset \Sc_i, \; i = 1,...,C} \prod_i p(\Zc_i|\yv) \left[ \begin{array}{c} \EE[\xv|\yv,\Zc_1] \\ \EE[\xv|\yv,\Zc_2] \\ \vdots \\ \EE[\xv|\yv,\Zc_C] \end{array} \right] \\
&=& \label{x_ammse2} \left[\begin{array}{c}
  \sum_{\Zc_1 \subset \Sc_1} p(\Zc_1|\yv) \EE[\xv|\yv,\Zc_1] \\ \sum_{\Zc_2 \subset \Sc_2} p(\Zc_2|\yv) \EE[\xv|\yv,\Zc_2] \\  \vdots \\
  \sum_{\Zc_C \subset \Sc_C} p(\Zc_C|\yv) \EE[\xv|\yv,\Zc_C]
\end{array} \right]
\end{eqnarray}
where the last line follows from the fact that
$\sum_{\Zc_i} p(\Zc_i|\yv) = 1$.
Thus, the semi-orthogonality of the columns in the sensing matrix allows us to obtain the MMSE estimate of $\xv$ in a divide-and-conquer manner by estimating the non-overlapping sections of $\xv$ independently from each other. Other structural properties of $\Psim$ can be utilized to reduce further the complexity of the MMSE estimation. For example, the orthogonal clusters exhibit some form of similarity and the columns within a particular cluster are also related to each other. We explore these properties for complexity reduction in Section \ref{sec:complexity}. However, before doing so, we devote the following section to a full description of our Bayesian orthogonal clustering algorithm.

\section{An Orthogonal Clustering (OC) Algorithm for Sparse Reconstruction}\label{orth_algo}
In this section, we present our sparse reconstruction algorithm, which is based on orthogonal clustering. The main steps of the algorithm are detailed in the following and summarized in~Figure~\ref{flowchart11}.

\subsection{Determine dominant positions}
Consider the model given in (\ref{model}) reproduced here for convenience,
$\yv = \Psim \xv + \nv.$
By correlating the observation vector, $\yv$, with the columns of the sensing matrix, $\Psim$, and by retaining correlations that exceed a certain threshold, we can determine the dominant positions/regions where the support of the sparse vector, $\xv$, is located. The performance of our orthogonal clustering algorithm is dependent on this initial correlation-based guess.\footnote{We can also apply a convex relaxation approach, retain the $P$ largest values, and form clusters around them. This allows us to incorporate {\em a priori} statistical information and obtain MMSE estimates but the algorithm in this case is bottle-necked by the performance of the convex relaxation approach and also loses the appeal of low complexity.}

\subsection{Form semi-orthogonal clusters}\label{form_clusters}

Define a threshold $\kappa$ such that $p(\nv>\kappa)$ $\define p_{\nv}$ is very small.\footnote{As $\nv \sim \Nc(0,\sigma_n^2)$, the threshold can be easily evaluated as, $\kappa = \sqrt{2\sigma_n^2}{\rm erfc}^{-1}(2p_\nv)$.} The previous correlation step creates a vector of $N$ correlations. From these correlations, obtain the indices with the correlation greater than the threshold, $\kappa$. Let $i_1$ denote the index with the largest correlation above $\kappa$ and form a cluster of size $L$ centered around $i_1$.\footnote{Given a fat sensing matrix, we consider two columns to be orthogonal (or semi orthogonal) when their correlation is below some value, $\varepsilon$. The cluster size $L$ is thus the minimum separation between two columns that makes these two columns semi-orthogonal. Obviously, the distance ,$L$, is a function of the correlation tolerance, $\varepsilon$. The lower the tolerance, $\varepsilon$, the larger the cluster size, $L$.}  Now, let $i_2$ denote the corresponding index of the second largest correlation above $\kappa$ and form another cluster of size $L$ around $i_2$. If the two clusters thus formed are overlapping, merge them into one big cluster. Continue this procedure until all the correlations greater than $\kappa$ are exhausted.


\subsection{Find the dominant supports and their likelihoods}
Let $L_i$ be the length of cluster $i$ and let $P_c$ denote the maximum possible support size in a cluster.\footnote{$P_c$ is calculated in a way similar to $P$ as the support in a cluster is also a Binomial distribution $\sim$ B($L_i,p$). Thus, we set $P_c = \lceil {\rm{erfc}}^{-1} (10^{-2}) \sqrt{2L_ip(1-p)} + L_ip \rceil$ (see footnote 5).} Let $C$ be the total number of semi-orthogonal clusters formed in the previous step. For each of them, find the most probable support of size, $|\Sc| = 1, |\Sc| = 2, \cdots, |\Sc| = P_c$, by calculating the likelihoods for all supports of size $|\Sc|$ (using either (\ref{map-metric}) or (\ref{map-metric1})). Each cluster is processed independently by capitalizing on the semi-orthogonality between the clusters. The expected value of the sparse vector $\xv$ given $\yv$ and the most probable support for each size can also be evaluated using either (\ref{Exp_MMSE}) or (\ref{Exp_LS}) depending on the {\em a priori} statistical~information.

\subsection{Evaluate the estimate of $\xv$}
Once we have the dominant supports for each cluster, their likelihoods, the expected value of $\xv$ given $\yv$ and the dominant supports, the MMSE (or MAP) estimates of $\hat{\xv}$ can be evaluated as discussed in Section \ref{sec:cont-band} (see (\ref{x_ammse2})). Note that these estimates are approximate as they are evaluated using only the dominant supports instead of using all supports.


\section{Reducing the Computational Complexity}\label{sec:complexity}

In this paper, we explore three structures of the sensing matrix that help us to reduce the complexity of MMSE estimation.
\begin{enumerate}
\item {\em Orthogonality (independence) of clusters}: In Section \ref{sec:cont-band}, the orthogonality of clusters allowed us to calculate the MMSE estimate independently over clusters in a divide-and-conquer manner.
\item {\em Similarity of clusters}: While the columns of the clusters are (semi)orthogonal, allowing us to treat them independently, these columns could exhibit some form of similarity making some MMSE calculations invariant over these clusters. For example, the columns of a DFT matrix can be obtained from each other through a modulation operation while those of the Toeplitz matrix can be obtained through a shift operation. The correlation calculations that repeatedly appear in the MMSE estimation are invariant to the modulation and shift operations.
\item {\em Order within a cluster}: MMSE estimation in a cluster involves calculating the likelihoods and expectations for all supports of size $i = 1, 2, \cdots, P_c$. Several quantities involved in these evaluations can be obtained in an order-recursive manner, incrementally moving from calculations for supports of size $i$ to similar calculations for supports of size $i+1$.
\end{enumerate}
We explore the last two properties in the following subsections.

\subsection{Similarity of Clusters} \label{from_one_cluster to other}
As evident from the previous sections, calculating the likelihood can be done in a divide-and-conquer manner by calculating the likelihood for each cluster independently. This is a direct consequence of the semi-orthogonality structure of the columns of the sensing matrix. 
Moreover, due to the rich structure of the sensing matrix, the clusters formed are quite similar. In the following subsections, we use the structure present in DFT and Toeplitz sensing matrices to show that the likelihood and expectation expressions in each cluster (for both the Gaussian and non-Gaussian cases) are strongly related, allowing many calculations across clusters to be shared.

\subsubsection{Discrete Fourier Transform (DFT) Matrices}
Let
$\psim_1,\; \psim_2, \; \cdots,\; \psim_{L}$ denote the sensing columns associated
with the first cluster. Then, it is easy to see that the corresponding columns for the $i^{th}$
cluster of equal length that are $\triangle_i$ columns away are,
$\psim_1 \odot \psim_{\triangle_i}, \; \psim_2 \odot \psim_{\triangle_i}, \; \cdots, \; \psim_{L} \odot \psim_{\triangle_i}$,
where $\psim_{\triangle_i}$ is some constant vector that depends on the sensing columns.\footnote{For example, if we use the last $M$ rows of the DFT matrix to construct the sensing matrix, then  
$\psim_{\triangle_i} = \left[ \exp{\left(-\frac{\jmath 2\pi (N-M)}{N} \triangle_i \right)} \; \exp{\left(-\frac{\jmath 2\pi (N-(M-1))}{N} \triangle_i \right)}\; \cdots \; \exp{\left(-\frac{\jmath 2\pi (N-1)}{N} \triangle_i \right)} \right]^\transp$.}
Assume that we evaluate the likelihood, $p(\Zc_1|\yv)$, and expectation, $\EE[\xv|\yv,\Zc_1]$, for a set of columns, $\Zc_1$, in the first cluster. For this set, we make the assumption that
\begin{equation}\label{red_com1}
\yv = \Psim_{\Zc_1} \xv + \nv.
\end{equation}
Now, let $\Zc_i$ denote the same set of columns chosen from the $i^{th}$ cluster that is $\triangle_i$ columns away (in other words $\Zc_i = \Zc_1 + \triangle_i$). For this set, we assume that
\begin{equation}\label{red_com2}
\yv = \Psim_{\Zc_i} \xv + \nv.
\end{equation}
Now (Hadamard) multiply 
both sides of the above equation by $\psim_{\triangle_i}^*$ to get
\begin{equation}\label{red_com3}
\yv \odot \psim_{\triangle_i}^* = \Psim_{\Zc_1} \xv + \nv \odot \psim_{\triangle_i}^*.
\end{equation}
Note that (\ref{red_com1}) and (\ref{red_com3}) have the same sensing matrix and the same noise statistics ($\nv$ is a white circularly symmetric Gaussian and hence is invariant to multiplication by $\psim_{\triangle_i}^*$). The only difference is that $\yv$ is modulated by the vector $\psim_{\triangle_i}^*$ in moving from the first to the $i^{th}$ cluster. This allows us to write
\begin{eqnarray}
p(\Zc_i|\yv) = p(\Zc_1|\yv \odot \psim_{\triangle_i}^*) \quad \mbox{and} \quad
\EE[\xv|\yv,\Zc_i] = \EE[\xv|\yv \odot \psim_{\triangle_i}^*, \Zc_1]
\end{eqnarray}
which is valid for both the Gaussian and non-Gaussian cases. In other words, if $\Zc_i$ is obtained from $\Zc_1$ by a constant shift, then any $\yv$-independent calculations remain the same while any calculations involving $\yv$ are obtained by modulating $\yv$ by the vector $\psim_{\triangle_i}^*$ as shown in Figure \ref{flowchart2}. For example, the likelihood in the Gaussian case reads
\begin{eqnarray}
p(\yv|\Zc_i) = \frac{\exp\left(-\|\yv\|^2_{\Sigmam_{\Zc_i}^{-1}} \right)}{\det(\Sigmam_{\Zc_i})}
= \frac{\exp\left(-\|\yv \odot \psim_{\triangle_i}^*\|^2_{\Sigmam_{\Zc_1}^{-1}} \right)}{\det(\Sigmam_{\Zc_1})}
\end{eqnarray}
and, in the non-Gaussian case, it reads
\begin{eqnarray}
p(\yv|\Zc_i) \simeq \exp \left( -\| \yv \|^2_{\Pm_{\Zc_i}^\perp}\right) = \exp \left( -\| \yv \odot \psim_{\triangle_i}^* \|^2_{\Pm_{\Zc_1}^\perp}\right).
\end{eqnarray}
We observe similar behavior in calculating the expectation. Thus, in the Gaussian case, we have
\begin{eqnarray}
\EE[\xv|\yv,\Zc_i] = \sigma_x^2 \Psim_{\Zc_i}^\herm \Sigmam_{\Zc_i}^{-1} \yv = \sigma_x^2 \Psim_{\Zc_1}^\herm \Sigmam_{\Zc_1}^{-1} (\yv \odot \psim_{\triangle_i}^*)
\end{eqnarray}
and in the non-Gaussian case, we have
\begin{eqnarray}
\EE[\xv|\yv,\Zc_i] = \left( \Psim_{\Zc_i}^\herm  \Psim_{\Zc_i} \right)^{-1} \Psim_{\Zc_i}^\herm \yv = \left( \Psim_{\Zc_1}^\herm  \Psim_{\Zc_1} \right)^{-1} \Psim_{\Zc_1}^\herm (\yv \odot \psim_{\triangle_i}^*).
\end{eqnarray}

\subsubsection{Toeplitz/Hankel Matrices}
In the Toeplitz or block Toeplitz case, the sensing matrix reads $\Psim = \left[ \Psim_{\Sc_1} \; \Psim_{\Sc_2} \; \cdots \; \Psim_{\Sc_C} \right]$.
Now, the clusters can be modified to make sure that they are identical (by stretching their end points if necessary) such that
\label{toep}
$\Psim_{\Sc_i} = [\begin{array}{ccccccc} \Om & \cdots & \Om & \Thetam^\transp & \Om & \cdots & \Om \end{array}]^\transp.$
In other words, the $\Psim_{\Sc_i}$s are simply shifted versions of each other. We now calculate the quantities $\det(\Sigmam_{\Zc_1})$, $\left\| \yv \right\|_{\Sigmam_{\Zc_1}^{-1}}^2$, and $\left\| \yv \right\|_{\Pm_{\Zc_1}^\perp}^2$ for a set $\Zc_1$ of columns of the first cluster. We then choose an identical set of columns, $\Zc_i$, in the $i^{th}$ cluster. Then, it is intuitively clear that
\begin{eqnarray}
\det(\Sigmam_{\Zc_i}) = \det(\Sigmam_{\Zc_1}), \quad
\left\| {\yv} \right\|_{\Sigmam_{\Zc_i}^{-1}}^2 = \left\| \yv \odot \wv_i \right\|_{\Sigmam_{\Zc_1}^{-1}}^2, \quad
\mbox{and} \quad \left\| {\yv} \right\|_{\Pm_{\Zc_i}^\perp}^2 = \left\| \yv \odot \wv_i \right\|_{\Pm_{\Zc_1}^\perp}^2
\end{eqnarray}
where $\wv_i$ is a rectangular window corresponding to the location of the non-zero rows of $\Psim_{\Sc_i}$.

\subsection{Order within a cluster} \label{Within_a_cluster}
To evaluate the likelihood for supports of size $i = 1, 2, ..., P_c$ in a single cluster, we pursue an order-recursive approach, calculating the likelihood and expectation for supports of size $i+1$ by updating calculations made for supports of size $i$. In the following, we assume that we have calculated the likelihood and expectation involving the columns, $\Psim_\Sc$, which we would like to update to $\Psim_{\Sc'} = [\Psim_\Sc \; \psim_i]$.

\subsubsection{$\xv|\Sc$ is Gaussian}
To calculate the likelihood
$\calL_{\Sc'} = \frac{\exp\left( -\frac{1}{\sigma_n^2} \| \yv \|^2_{\Sigmam_{\Sc'}^{-1}} \right)}{\det(\Sigmam_{\Sc'})}$
with $\Sigmam_{\Sc'} = \Id_M + \frac{\sigma_x^2}{\sigma_n^2} \Psim_{\Sc'} \Psim_{\Sc'}^\herm$, note that
$\Sigmam_{\Sc'} = \Sigmam_\Sc + \frac{\sigma_x^2}{\sigma_n^2} \psim_{i} \psim_{i}^\herm,$
or by the matrix inversion lemma,
\begin{eqnarray}
\label{mil_Sigma} \Sigmam^{-1}_{\Sc'} &=& \Sigmam^{-1}_{\Sc} - \frac{\sigma_x^2}{\sigma_n^2} \xi_i \omegam_i \omegam_i^\herm
\end{eqnarray}
where
\begin{eqnarray}
\label{omega_G} \omegam_i &\define& \Sigmam_{\Sc}^{-1} \psim_{i}\\
\label{xi_G} \xi_i &\define& \left(1 + \frac{\sigma_x^2}{\sigma_n^2} \psim_{i} ^\herm \Sigmam_{\Sc}^{-1} \psim_{i} \right)^{-1} =  \left(1 + \frac{\sigma_x^2}{\sigma_n^2} \psim_{i} ^\herm \omegam_i \right)^{-1}.
\end{eqnarray}
As we are actually interested in computing $\exp\left( -\frac{1}{\sigma_n^2} \| \yv \|^2_{{\Sigmam^{-1}_{\Sc'}}} \right)$, using (\ref{mil_Sigma}) we obtain
\begin{eqnarray}
\nonumber \exp\left( -\frac{1}{\sigma_n^2} \| \yv \|^2_{{\Sigmam^{-1}_{\Sc'}}} \right) &=& \exp\left( -\frac{1}{\sigma_n^2} \| \yv \|^2_{{\Sigmam^{-1}_{\Sc}}}
+  \frac{\sigma_x^2\xi_i}{\sigma_n^4} \|\omegam_i^\herm \yv \|^2 \right) \\
&=& \label{inv_Sigma} \exp\left( -\frac{1}{\sigma_n^2} \| \yv \|^2_{{\Sigmam^{-1}_{\Sc}}} \right) \exp \left( \frac{\sigma_x^2\xi_i}{\sigma_n^4} \|\omegam_i^\herm \yv \|^2 \right).
\end{eqnarray}
The determinant of $\Sigmam_{\Sc'}$ can be evaluated 
as follows:
\begin{eqnarray}
\label{det_Sigma} \det(\Sigmam_{\Sc'}) = \det \left(\Sigmam_{\Sc} + \frac{\sigma_x^2}{\sigma_n^2} \psim_{i} \psim_{i}^\herm\right)
= \det \left(1 + \frac{\sigma_x^2}{\sigma_n^2} \psim_{i} ^\herm \Sigmam_{\Sc}^{-1} \psim_{i} \right) \; \det \left(\Sigmam_{\Sc}\right)
=  \xi_i^{-1} \; \det \left(\Sigmam_{\Sc}\right).
\end{eqnarray}
Thus, the likelihood for the support of size $\Sc'$ can be written as (using (\ref{inv_Sigma}) and (\ref{det_Sigma})),
\begin{eqnarray}
\label{like_G} \nonumber \calL_{\Sc'} &=& \frac{\exp\left( -\frac{1}{\sigma_n^2} \| \yv \|^2_{{\Sigmam^{-1}_{\Sc}}} \right) \exp \left(  \frac{\sigma_x^2\xi_i}{\sigma_n^4} \|\omegam_i^\herm \yv \|^2 \right)}{\det \left(\Sigmam_{\Sc}\right)\xi_i^{-1}} \\
 &=& \calL_{\Sc} \;\; {\underbrace{ \xi_i\;{\exp \left(  \frac{\sigma_x^2\xi_i}{\sigma_n^4} \|\omegam_i^\herm \yv \|^2 \right)}}}.\\
&& \nonumber \quad \quad \quad \quad \quad \quad \; \delta_{i}
\end{eqnarray}
This shows that to calculate $\calL_{\Sc'}$, we need to compute only $\omegam_i$ and $\xi_i$, which constitute $\delta_{i}$. To calculate $\omegam_i$ for a cluster of length, $L$, $\calO(LM^2)$ operations is required if standard matrix multiplication is used. This complexity can be reduced to $\calO(LM)$ by storing all the past computed values of $\omegam$ and $\xi$ and using the structure of $\Sigmam_\Sc$ \cite{FBMP}.

Similarly, $\EE[\xv_{\Sc'}|\yv]$ can be calculated in an order-recursive manner as follows:
\begin{eqnarray}\label{Exp_i_G}
\EE[\xv_{\Sc'}|\yv] = \left[ \begin{array}{c}
    \EE[\xv_\Sc|\yv] \\
  \sigma_x^2 \omegam_i^\herm \yv
\end{array} \right].
\end{eqnarray}



\subsubsection{$\xv|\Sc$ is unknown}

To calculate the likelihood in the non-Gaussian case, we need to evaluate the norm,
$\left\| \yv \right\|^2_{\Pm_{\Sc'}^\perp}=\left\| \yv \right\|^2 - \yv^\herm \Psim_{\Sc'} \left( \Psim_{\Sc'}^\herm \Psim_{\Sc'} \right)^{-1} \Psim_{\Sc'}^\herm \yv.$
Our approach mainly hinges on calculating the inverse $\Lambdam_{\Sc'} \define \left( \Psim_{\Sc'}^\herm \Psim_{\Sc'} \right)^{-1}$ recursively. We do this by invoking the block inversion formula
\begin{eqnarray}
\label{red_un1}
\Lambdam_{\Sc'} &=& \left[ \begin{array}{cc}
\Lambdam_{\Sc} + \frac{1}{\xi_i}\omegam_i \omegam_i^\herm  & -\frac{1}{\xi_i}\omegam_i \\
-\frac{1}{\xi_i}\omegam_i^\herm & \frac{1}{\xi_i}
\end{array}\right]\\
\end{eqnarray}
where
\begin{eqnarray}
\label{red_un2} 
\omegam_i &\define& \Lambdam_{\Sc} (\Psim_\Sc^\herm \psim_i) \\
\label{red_un3} \xi_i &\define& \left\| \psim_i \right\|^2 - (\psim_i^\herm \Psim_\Sc) \Lambdam_{\Sc} (\Psim_\Sc^\herm \psim_i) = \left\| \psim_i \right\|^2 - \omegam_i^\herm \etam_i
\end{eqnarray}
with the elements of $\etam_i \define \Psim_\Sc^\herm \psim_i$ all available (i.e., they are calculated initially and can be reused afterwards). Using this recursion, we can construct (following some straightforward manipulation) a recursion for the projected norm $\calL_{\Sc'}$:
\begin{eqnarray}\label{red_un4}
\nonumber \calL_{\Sc'} &=& \exp\left( -\frac{1}{\sigma_n^2} \left[ \left\| \yv \right\|^2 - \yv^\herm \Psim_{\Sc'} \Lambdam_{\Sc'} \Psim_{\Sc'}^\herm \yv \right] \right)\\
&=& \nonumber \exp\left( -\frac{1}{\sigma_n^2} \left[ \left\| \yv \right\|^2 - \yv^\herm \Psim_{\Sc} \Lambdam_{\Sc} \Psim_{\Sc}^\herm \yv \right]\right) \\
&& \nonumber \exp\left(-\frac{1}{\sigma_n^2} \left[
-\frac{1}{\xi_i} | (\yv^\herm \Psim_\Sc) \omegam_i|^2 + \frac{2}{\xi_i} \Re\{ (\yv^\herm \psim_i) \omegam_i^\herm (\Psim_\Sc^\herm \yv)\} - \frac{1}{\xi_i} |\yv^\herm \psim_i|^2  \right] \right) \\
&=& \calL_\Sc \; \underbrace{ \exp\left(\frac{1}{\sigma_n^2 \xi_i} \left[
| (\yv^\herm \Psim_\Sc) \omegam_i|^2 - 2 \Re\{ (\yv^\herm \psim_i) \omegam_i^\herm (\Psim_\Sc^\herm \yv)\} + |\yv^\herm \psim_i|^2  \right] \right)}.\\
&& \nonumber \quad \quad \quad \quad \quad \quad \quad \quad \quad \quad \quad\quad \quad \quad \quad\; \delta_{i}
\end{eqnarray}

Similarly, we can show that
\begin{eqnarray}
  \label{Exp_i_NG}  \EE[\xv_{\Sc'}|\yv] = \Lambdam_{\Sc'}(\Psim_{\Sc'}^\herm \yv) =
  \left[ \begin{array}{c}
    \EE[\xv_{\Sc}|\yv] + \frac{1}{\xi_i}\omegam_i \etam_i^\herm \EE[\xv_{\Sc}|\yv] - \frac{1}{\xi_i}\omegam_i\psim_i^\herm \yv \\
    -\frac{1}{\xi_i} \etam_i^\herm \EE[\xv_{\Sc}|\yv] + \frac{1}{\xi_i} \psim_i^\herm \yv
  \end{array} \right].
\end{eqnarray}

The cluster independent and cluster-wise evaluations in our recursive procedure for both the cases ($\xv|\Sc$ Gaussian or unknown) are summarized in Table \ref{t2}.

\begin{table}[h]
\begin{center}
\caption{Cluster independent and cluster-wise evaluations involved in the recursive procedure for complexity reduction within a cluster}
\begin{tabular}{|l||l|l|}
\hline & {\small \bf Cluster Independent Evaluations} & {\small \bf Cluster-wise Evaluations} \\
\hline
\hline & \small Evaluate $\omegam_i$ and $\xi_i$ using (\ref{omega_G}) and (\ref{xi_G}) & \small Evaluate $\| \omegam_i^\herm \yv \|^2$ \\
\small $\xv|\Sc$ is Gaussian & \small Update $\Sigmam_{\Sc}^{-1}$ using (\ref{mil_Sigma}) & \small Update $\calL_{\Sc}$ using equation (\ref{like_G}) \\
& \small Update $\det(\Sigmam_{\Sc})$ using (\ref{det_Sigma}) & \small Update $\EE[\xv_{\Sc'}|\yv]$ using equation (\ref{Exp_i_G})\\
\hline & \small Initialize: Calculate $\psim_i^\herm \psim_j \; \forall \; i,j$ & \small Initialize: Evaluate $\yv^\herm \psim_i \; \forall \; i$ \\
\small $\xv|\Sc$ is unknown & \small Evaluate $\omegam_i$ using equation (\ref{red_un2}) & \small Update $\calL_{\Sc}$ using equation (\ref{red_un4}) \\
& \small Update $\Sigmam_{\Sc}$ using equations (\ref{red_un1}) and (\ref{red_un3}) & \small Update $\EE[\xv_{\Sc'}|\yv]$ using equation (\ref{Exp_i_NG}) \\
\hline
\end{tabular}\label{t2}
\end{center}
\end{table}

\section{Simulation Results}\label{sec:sim}

In this section, we compare the performance of the OC algorithm with popular sparse reconstruction methods available in the literature including the convex relaxation (CR) method \cite{candes-romberg-tao}, OMP \cite{OMP2}, and FBMP \cite{FBMP}. The parameters of these algorithms are set according to the specifications provided by the authors to achieve the best results.\footnote{For a fair comparison, we perform the MMSE refinement on the output of CR and OMP.}
The parameters that we use in all the simulations are $N = 800$, $M = \frac{N}{4} = 200$, $p = 10^{-2}$, and SNR $= 30$dB (unless stated otherwise). Specifically, we demonstrate the performance of our algorithm for the case when the sensing matrix is a DFT or a Toeplitz matrix. We start by first investigating the effect of cluster length on the performance of OC.

\subsection{The effect of the cluster length, $L$}
Figure \ref{var_L} compares the normalized mean-square error (NMSE) of OC as the cluster length, $L$, is varied. The NMSE is defined as
${\rm NMSE} = \frac{1}{R} \sum_{r = 1}^{R} \frac{ \left\| \hat{\xv}^{(r)} - {\xv}^{(r)} \right\|^2}{ \left\| {\xv}^{(r)} \right\|^2}$,
where $\hat{\xv}$ stands for the estimated sparse signal for realization $r$, and $R$ is the total number of runs. For this case, the DFT matrix is used as the sensing matrix with $\xv|\Sc$ Gaussian. Note that while implementing OC with fixed-length clusters, overlapping of clusters is not allowed to maintain orthogonality. This results in an increase in the probability of missing the correct support if two supports are close to each other. Thus, the smaller the cluster, the greater the probability of missing the correct supports. This is evident from Figure \ref{var_L} as performance of OC improves by increasing $L$. Obviously, this improvement in performance is obtained at the expense of speed. Figure \ref{runtime_L} shows that the smaller the length of clusters, the faster the algorithm.
Note that for larger values of $L$ (e.g., $L>32$), it might not be possible to form the required number of non-overlapping clusters. To overcome this problem, we present the performance of OC implemented with variable length clusters (labeled as ``OC'' in Figure \ref{var_L}). In this case, the overlapping clusters are joined together to form larger clusters. It can be observed from Figure \ref{var_L} that the performance of OC with variable-length clusters is better than the case when it is implemented with fixed-length clusters. Moreover, this performance is achieved with a reasonable run-time{\footnote{Thus, the following simulation results are presented with OC implemented using variable length clusters.} as shown in Figure \ref{runtime_L}.

\subsection{The effect of the signal-to-noise ratio (SNR)}
Figure \ref{nmse_DFT_Gaussian} compares the performance of the algorithms for the case when the sensing matrix is a DFT matrix and $\xv|\Sc$ is Gaussian. 
In the FBMP implementation, the number of greedy branches to explore ($D$) is set to $10$. Note that OC outperforms all other algorithms at low SNR while FBMP performs quite close to it at SNR $\geq 25$ dB. It outperforms both OMP and CR at all SNR values. Specifically, at SNR $= 25$ dB, OC has a gain of approximately 2 dB and 3 dB over CR and OMP, respectively.
The performance of the algorithms for the case when the sensing matrix is a DFT matrix and $\xv|\Sc$ is unknown is presented in Figure \ref{nmse_DFT_nonGaussian}. In this case, the entries of $\xv_G$ are drawn from a uniform distribution. Here, FBMP is allowed to estimate the hyper-parameters using its approximate ML algorithm (with $E$ set to $10$)\cite{FBMP}. It can be seen that OC easily outperforms OMP and FBMP while CR performs similar to OC. Specifically, at SNR $= 25$ dB, OC outperforms OMP and FBMP by approximately $5$ dB. 
Figure \ref{nmse_Toeplitz_Gaussian} compares the performance of the algorithms for the case when the sensing matrix is Toeplitz. To do so, we first generate a Toeplitz matrix from a column having $20$ non-zero consecutive samples drawn from a Gaussian distribution. The sensing matrix is then extracted by uniformly sub-sampling this full matrix at a rate less than the duration of the signal.\footnote{In this case, the sub-sampling rate is 4 times less making $M = 200$.}
Note that the performance of OC and FBMP is almost the same at low SNR but OC outperforms FBMP in the high SNR region.
OMP and CR do not perform well in this case as the sensing matrix does not exhibit the requisite incoherence conditions (in this case, $\mu(\Psim) \simeq 0.9$) on which much of the CS theory is based.

\subsection{The effect of the under-sampling ratio $(\frac{N}{M})$}
Figure \ref{nmse_us_DFT_Gaussian} shows the performance of the algorithms (for the case when the sensing matrix is DFT and $\xv|\Sc$ is Gaussian) when the under-sampling ratio $(\frac{N}{M})$ is varied. It can be observed that the performance of all the algorithms deteriorates as $\frac{N}{M}$ increases. OC and FBMP perform quite close to each other with OC performing slightly better at high $(\frac{N}{M})$ ratios.

\subsection{The effect of the sparsity rate, $p$}
Figure \ref{nmse_p_DFT_Gaussian} compares the performance of the algorithms when the sparsity rate, $p$, is varied (for the case when the sensing matrix is DFT and $\xv|\Sc$ is Gaussian). It can be seen that the performance of OC is quite close to CR and FBMP at low sparsity rate while it outperforms OMP by approximately $3$ dB for the studied range of $p$. The performance of OC deteriorates at the high sparsity rate because the number of clusters increases as $p$ increases and the probability of clusters to be near or overlapping each other increases. Thus, in this case, the orthogonality assumption of OC becomes weak. Figure \ref{runtime_p_DFT_Gaussian} compares the mean run-time of all the algorithms. It can be seen that OC is faster than all other algorithms. 
As sparsity rate increases, the length of the clusters increases, and thus the complexity of OC. 
Figure \ref{nmse_p_DFT_nonGaussian} shows that OC performs quite well at the low sparsity rate in the case when the sensing matrix is DFT and $\xv|\Sc$ is unknown. FBMP does not perform well at the low sparsity rate in this case even with its approximate ML algorithm. The run-time of FBMP is also higher as compared to Figure \ref{nmse_p_DFT_Gaussian} due to the time taken to estimate the hyper-parameters using the ML algorithm.
In the case of the Toeplitz matrix (see Figure \ref{nmse_p_Toeplitz_Gaussian}), the performance of OC and FBMP is almost the same while the performance of CR and OMP is quite poor due to the weak incoherence of the sensing matrix. It can also be observed from Figure \ref{runtime_p_Toeplitz_Gaussian} that OC is quite fast compared to the other algorithms. 

\section{Conclusion and Future Work}\label{sec:conc}
In this paper, we present the Orthogonal Clustering algorithm for fast Bayesian sparse reconstruction. This algorithm makes collective use of the underlying structure (sparsity, {\em a priori} statistical information, structure of the sensing matrix) to achieve superior performance at much lower complexity compared with other algorithms especially at low sparsity rates. The proposed algorithm has the following distinctive features.
\begin{enumerate}
\item It is able to deal with Gaussian priors as well as with priors that are non-Gaussian or unknown.
\item It utilizes the structure of the sensing matrix, including orthogonality, modularity, and order-recursive calculations.
\item In the Gaussian case, OC beats all other algorithms in terms of complexity and performance for low sparsity rates. In the non-Gaussian case, it outperforms all other algorithms (most notably FBMP) for both low and high sparsity rates. Hence, the only disadvantage of OC is its performance at high sparsity rates. In this case, the clusters are no longer orthogonal, which results in large clusters and the orthogonality assumption becomes invalid. Fortunately, this drawback is only observed in the Gaussian case while in the non-Gaussian case, OC maintains a relative advantage over the other algorithms for all sparsity rates.
\item It is able to provide computable measures of performance (See \cite{impulse-noise} for details on how to calculate the error covariance matrix using orthogonality). 
\end{enumerate}
Our future work includes
\begin{enumerate}
\item The OC algorithm assumes that various clusters do not interact. We guarantee this by lumping any two clusters that are too close into a single larger cluster. This prevents us from implementing a fixed-size cluster algorithm and gives our algorithm the advantage of being computationally cleaner and more efficient. A prerequisite to do so however is to implement an OC that takes into account the interaction between neighboring  clusters.
\item The OC algorithm utilizes various levels of structure in the sensing matrix but falls short of utilizing one additional structure. Specifically, the various columns of any cluster are not random but are actually related (e.g., adjacent columns in the Toeplitz case exhibit a shift structure).\footnote{This structure is for example used in the lattice implementation of recursive least squares for drastic reduction in complexity \cite{Sayed}.} This additional structure can be used to reduce further the complexity of our algorithm.
\item The OC algorithm does not use any dependence between the active sparse elements (e.g., block sparsity). It can be specialized to deal with such cases.
\item The divide-and-conquer approach that we are able to pursue due to the structure of the sensing matrix can be utilized in the existing algorithms like OMP.
\end{enumerate}


\begin{figure}[h]
\centering
\subfloat[Normal]{
\includegraphics[width=.4\textwidth]{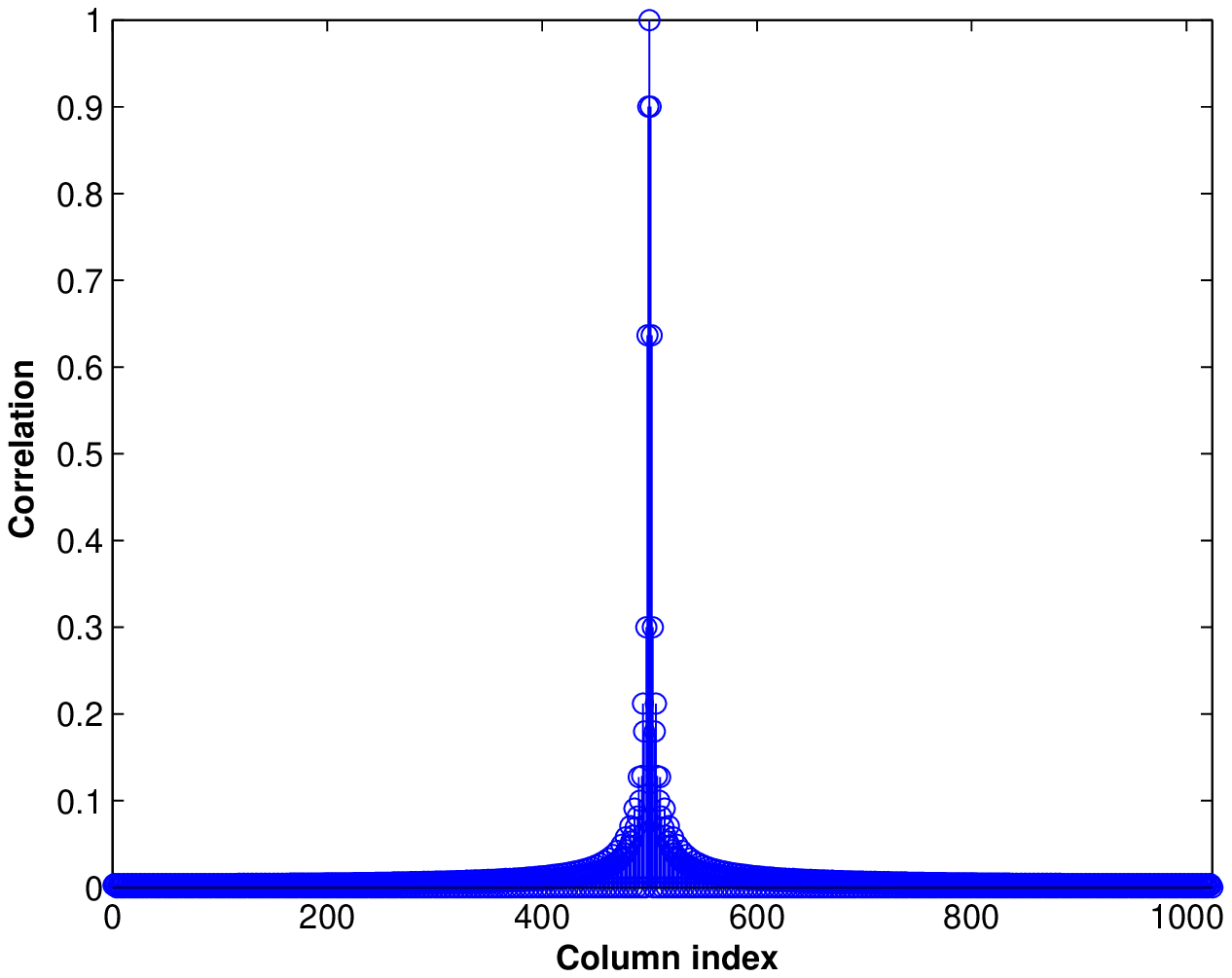}
\label{corr_N1024_DFT}
}
\subfloat[Zoomed]{
\includegraphics[width=.4\textwidth]{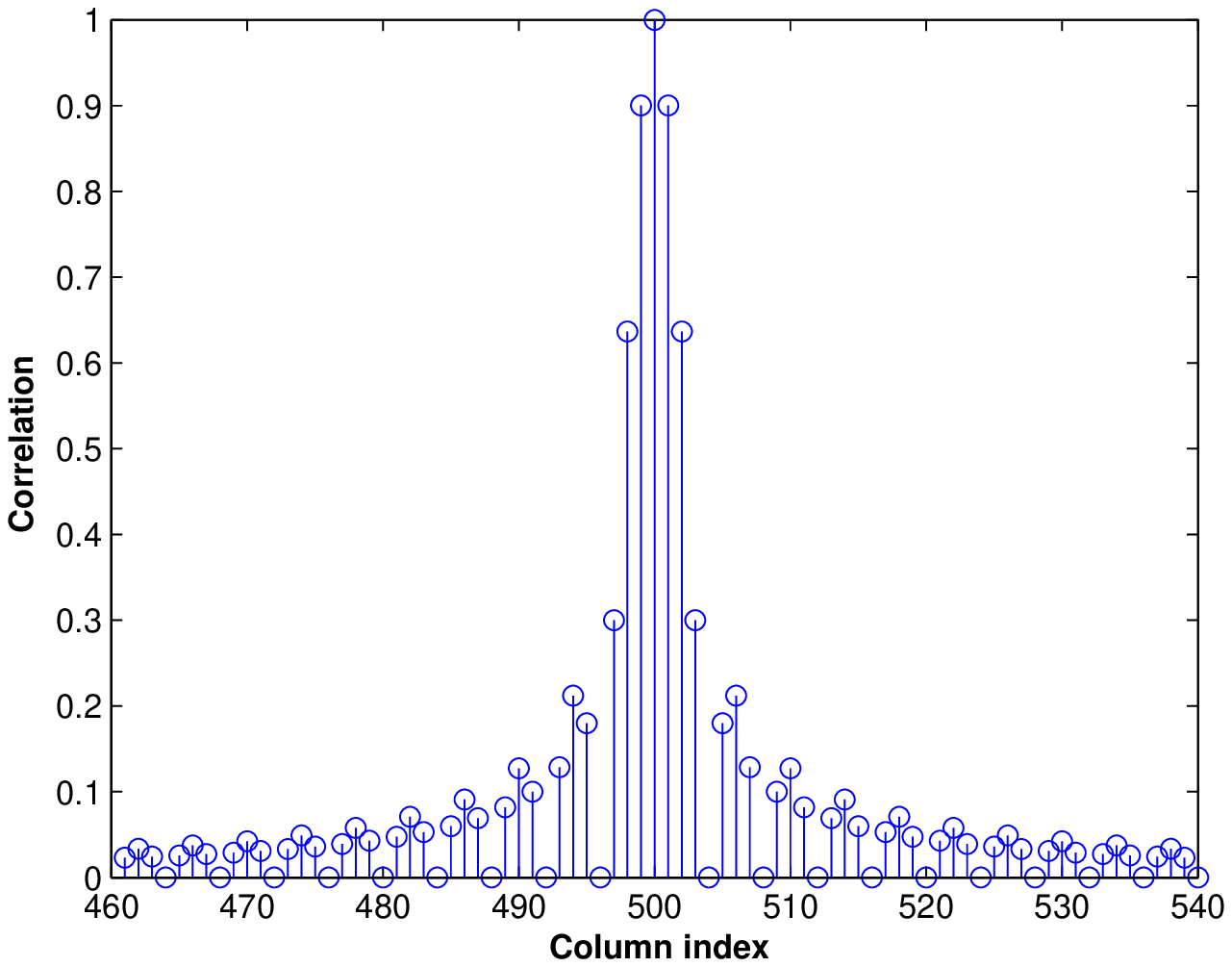}
\label{corr_N1024_DFT_zoomed}
}
\caption{\small The $500^{th}$ column has high correlation with its neighbors}
\label{corr_N_DFT}
\end{figure}

\begin{figure}[h]
\begin{minipage}[htp]{.5\textwidth}
{\begin{center}
\small
\psframebox[linearc=0.5,cornersize=absolute,framesep=10pt]{%
  \psset{shadowcolor=black!70,blur=true}%
  \begin{psmatrix}[rowsep=0.4,colsep=0.5]
    \psovalbox[shadow=true]{\parbox{1cm} {\centering \scriptsize Begin}} \\
    \psframebox[shadow=true]{\parbox{6cm} {\centering \scriptsize Correlate observation vector, $\yv$, with the sensing matrix, $\Psim$}}\\
    \psframebox[shadow=true]{\parbox{6cm} {\centering \scriptsize Form semi-orthogonal clusters around the positions with correlation values greater than the threshold, $\kappa$}}\\
    \psframebox[shadow=true]{\parbox{6cm} {\centering \scriptsize Process each cluster independently and in each cluster, calculate the
    likelihoods for supports of size, $|\Sc| = 1, |\Sc| = 2, \cdots, |\Sc| = P_c$}} \\
    \psframebox[shadow=true]{\parbox{6cm} {\centering \scriptsize Find the dominant supports of size, $|\Sc| = 1, |\Sc| = 2, \cdots, |\Sc| = P_c$, for each cluster}}\\
    \psframebox[shadow=true]{\parbox{6cm} {\centering \scriptsize Find $\EE[\xv|\yv,\Sc]$ for dominant support of each size}}\\
    \psframebox[shadow=true]{\parbox{6cm} {\centering \scriptsize Evaluate $\hat{\xv}_{\rm MMSE}$ or $\hat{\xv}_{\rm MAP}$}}\\
    \psovalbox[shadow=true]{\parbox{1cm} {\centering \scriptsize End}}
    \psset{linewidth=1.5pt}
    \ncline{->}{1,1}{2,1}
    \ncline{->}{2,1}{3,1}
    \ncline{->}{3,1}{4,1}
    \ncline{->}{4,1}{5,1}
    \ncline{->}{5,1}{6,1}
    \ncline{->}{6,1}{7,1}
    \ncline{->}{7,1}{8,1}
  \end{psmatrix}
}
\end{center}}
\caption{\small Flowchart of the OC algorithm}
\label{flowchart11}
\end{minipage}
\begin{minipage}[htp]{.5\textwidth}
{\begin{center}
\small
{%
  \psset{shadowcolor=black!70,blur=true}%
  \begin{psmatrix}[rowsep=0.6,colsep=1]
    ~ & \scriptsize $\yv$ \\
    ~ & \psframebox[shadow=true]{\parbox{2cm}{\centering \scriptsize Modulator \\ \scriptsize $\yv \odot \psim^{*}_{\triangle_i}$}} \\
    ~ & \psframebox[shadow=true]{\parbox{3.5cm}{\centering \scriptsize Calculate \\ \scriptsize $p(\Zc_i|\yv)$ and $\EE[\xv|\yv,\Zc_i]$}}
    \psset{linewidth=1pt}
    \ncline{->}{1,2}{2,2}
    \ncline{->}{2,1}{2,2}<{\scriptsize $\psim_{\triangle_i}$}
    \ncline{->}{2,2}{3,2}
    \ncline{->}{3,1}{3,2}^{{\hspace*{-4.5cm} \scriptsize $\det\left(\Sigmam_{\Zc_1}\right)$ \& $\Sigmam_{\Zc_1}^{-1}$}}<{\hspace*{4cm} \scriptsize (or ${\Pm}_ {\Zc_1}^\perp$)}
    \ncline{->}{3,1}{3,2}
\end{psmatrix}%
}
\end{center}}
\caption{\small Block diagram of the reduced complexity algorithm for the DFT matrix}
\label{flowchart2}
\end{minipage}
\end{figure}



\begin{figure}[t]
\begin{minipage}[htp]{.5\textwidth}
\begin{center}
\includegraphics[width=\textwidth]{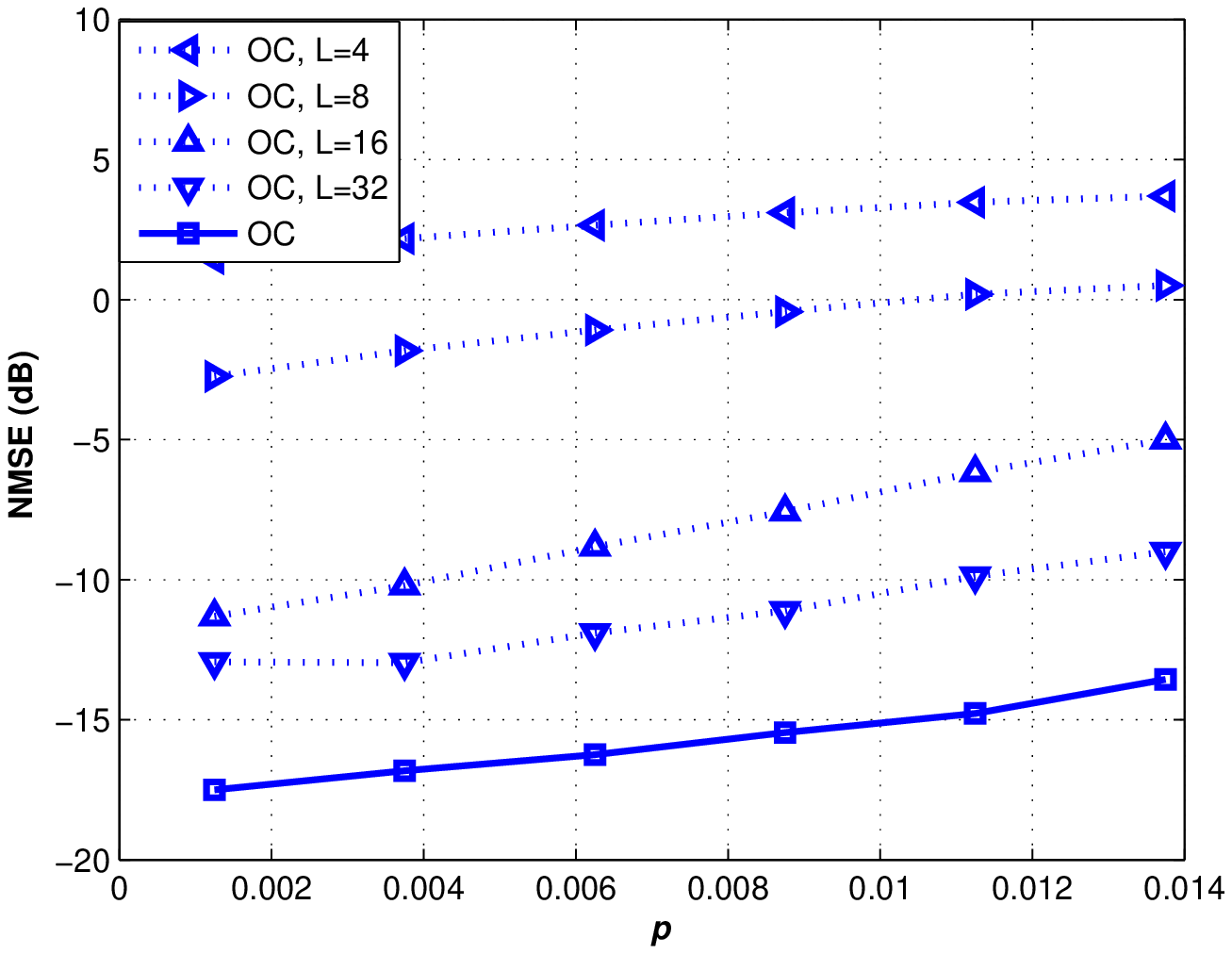}
\end{center}
\caption{\small NMSE vs $p$ for the OC algorithm with the length of the cluster varied.} \label{var_L}
\end{minipage}
\begin{minipage}[htp]{.5\textwidth}
\begin{center}
\includegraphics[width=\textwidth]{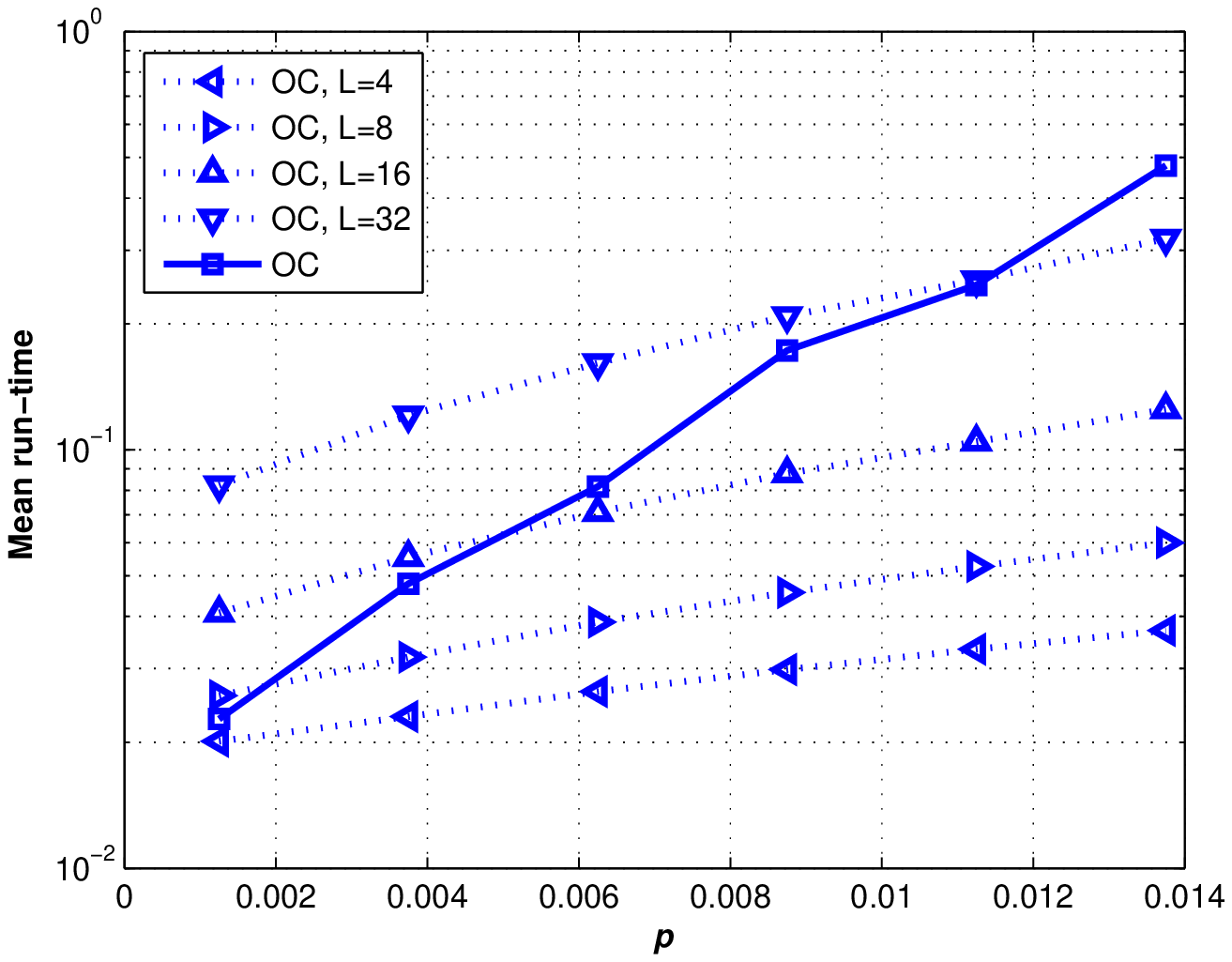}
\end{center}
\caption{\small Mean run-time for the OC algorithm with the length of cluster varied.} \label{runtime_L}
\end{minipage}
\begin{minipage}[htp]{.5\textwidth}
\begin{center}
\includegraphics[width=\textwidth]{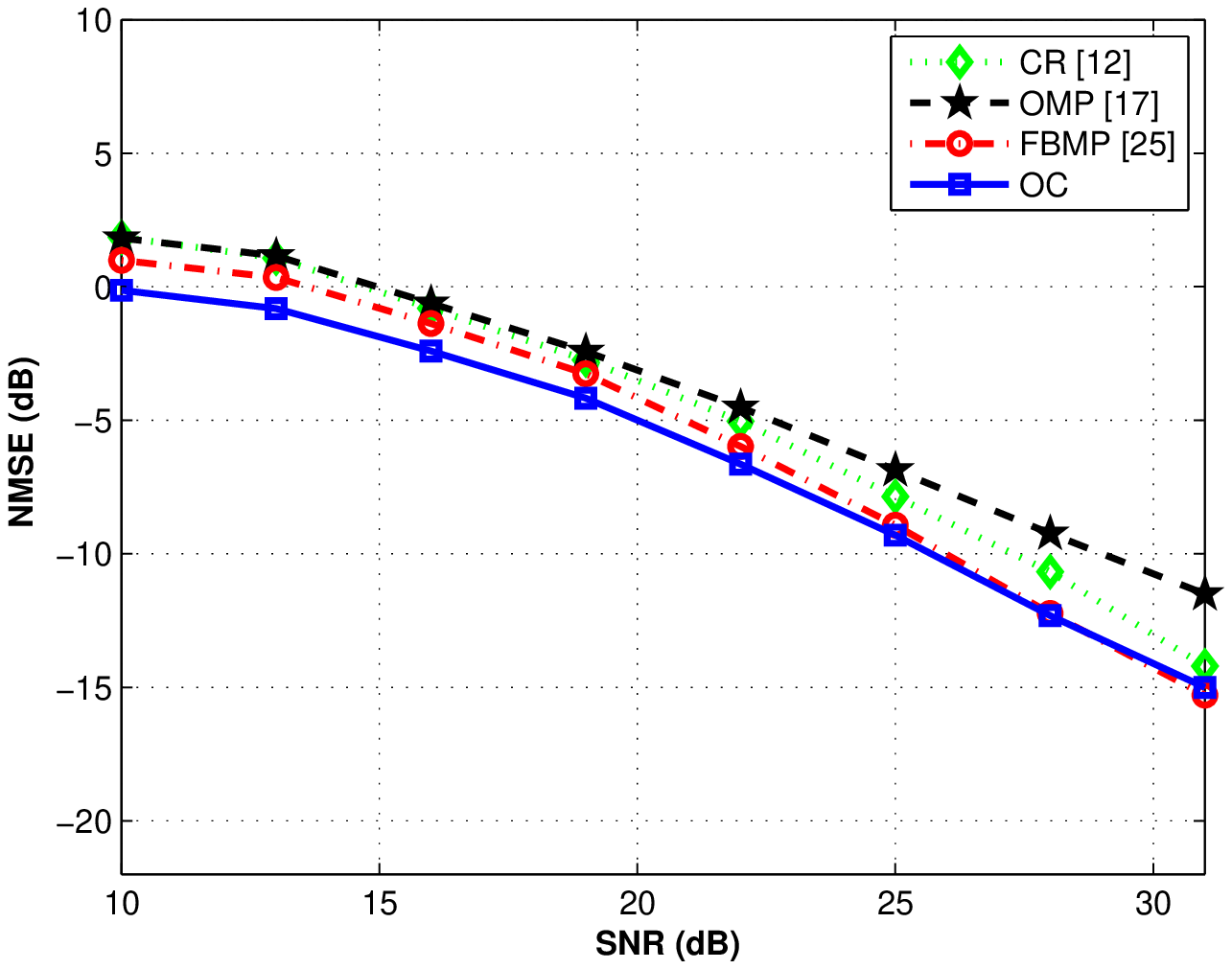}
\end{center}
\caption{\small NMSE vs SNR for the DFT matrix and $\xv|\Sc$ Gaussian.} \label{nmse_DFT_Gaussian}
\end{minipage}
\begin{minipage}[htp]{.5\textwidth}
\begin{center}
\includegraphics[width=\textwidth]{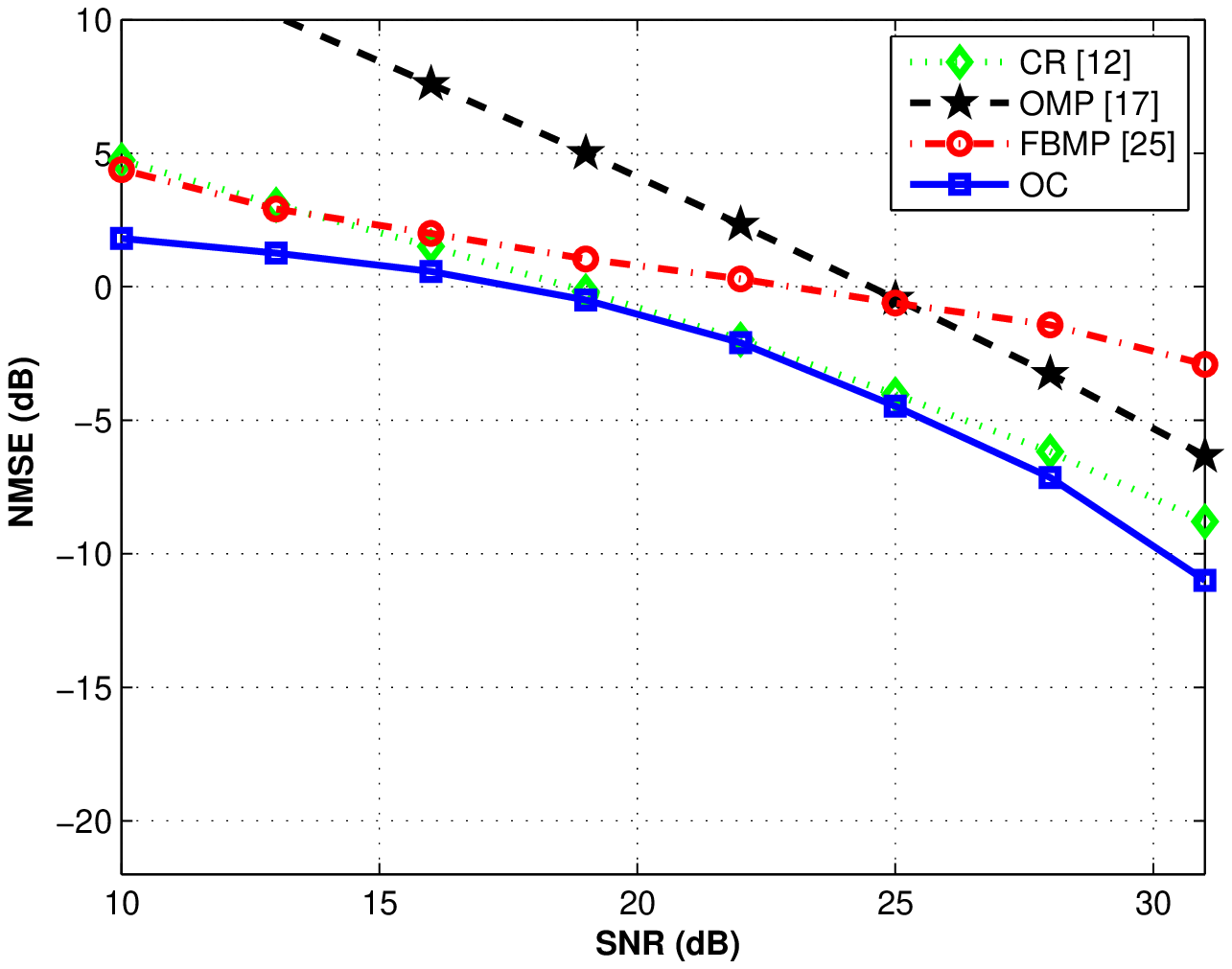}
\end{center}
\caption{\small NMSE vs SNR for the DFT matrix and $\xv|\Sc$ unknown.} \label{nmse_DFT_nonGaussian}
\end{minipage}
\begin{minipage}[htp]{.5\textwidth}
\begin{center}
\includegraphics[width=\textwidth]{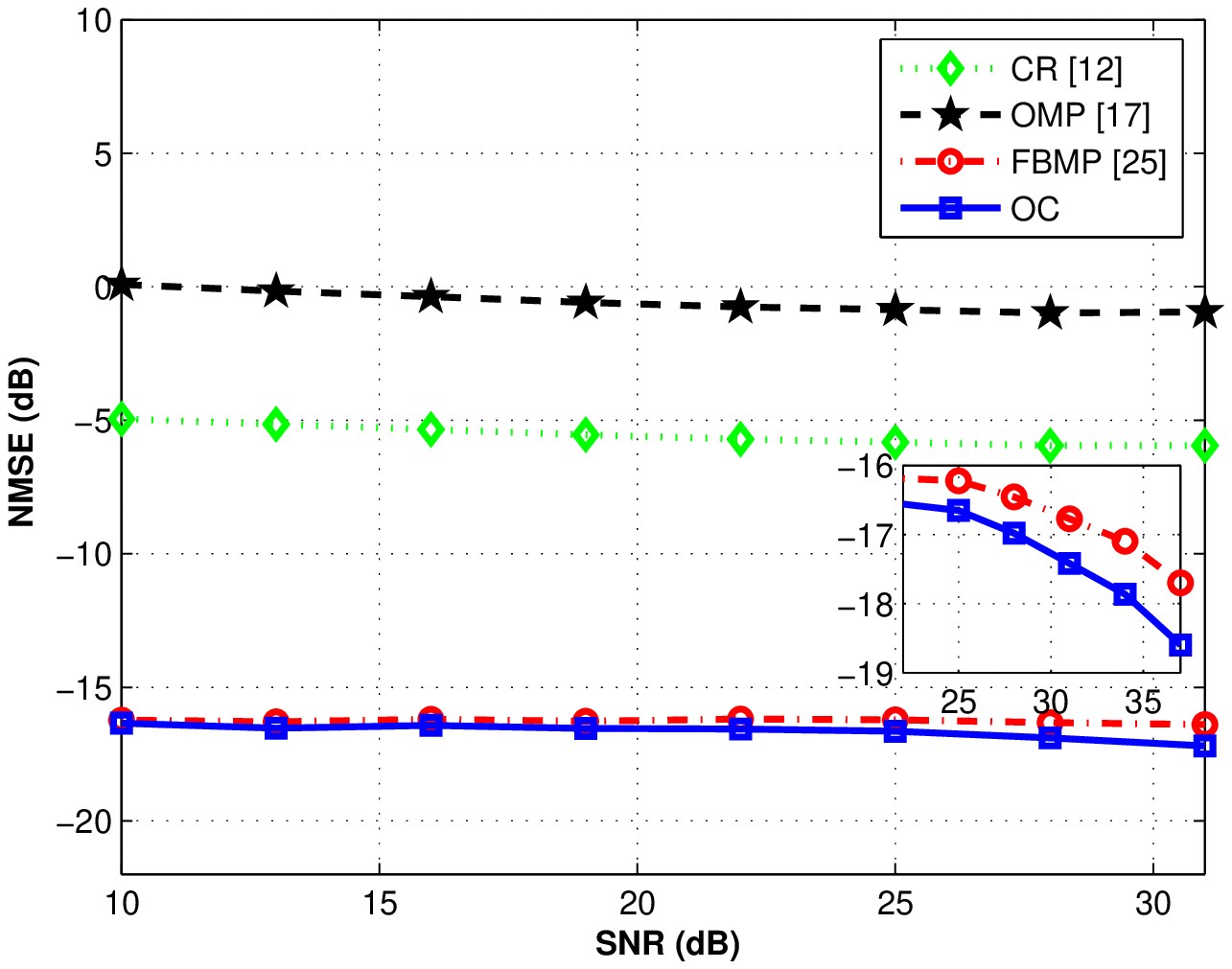}
\end{center}
\caption{\small NMSE vs SNR for the Toeplitz matrix and $\xv|\Sc$ Gaussian.} \label{nmse_Toeplitz_Gaussian}
\end{minipage}
\begin{minipage}[htp]{.5\textwidth}
\begin{center}
\includegraphics[width=\textwidth]{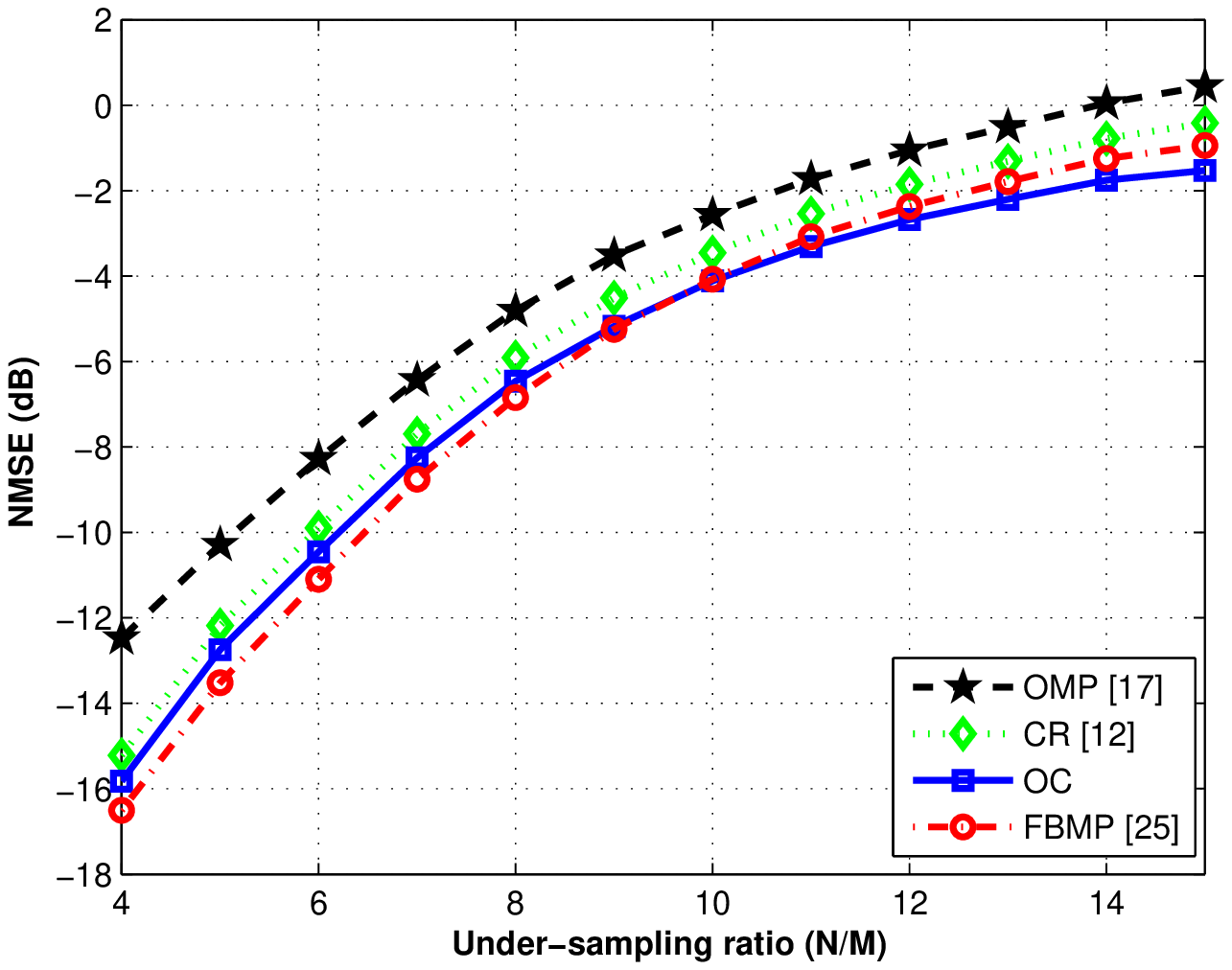}
\end{center}
\caption{\small NMSE vs the undersampling ratio ($\frac{N}{M}$) for the DFT matrix and $\xv|\Sc$ Gaussian.} \label{nmse_us_DFT_Gaussian}
\end{minipage}

\end{figure}

\begin{figure}[t]
\begin{minipage}[htp]{.5\textwidth}
\begin{center}
\includegraphics[width=\textwidth]{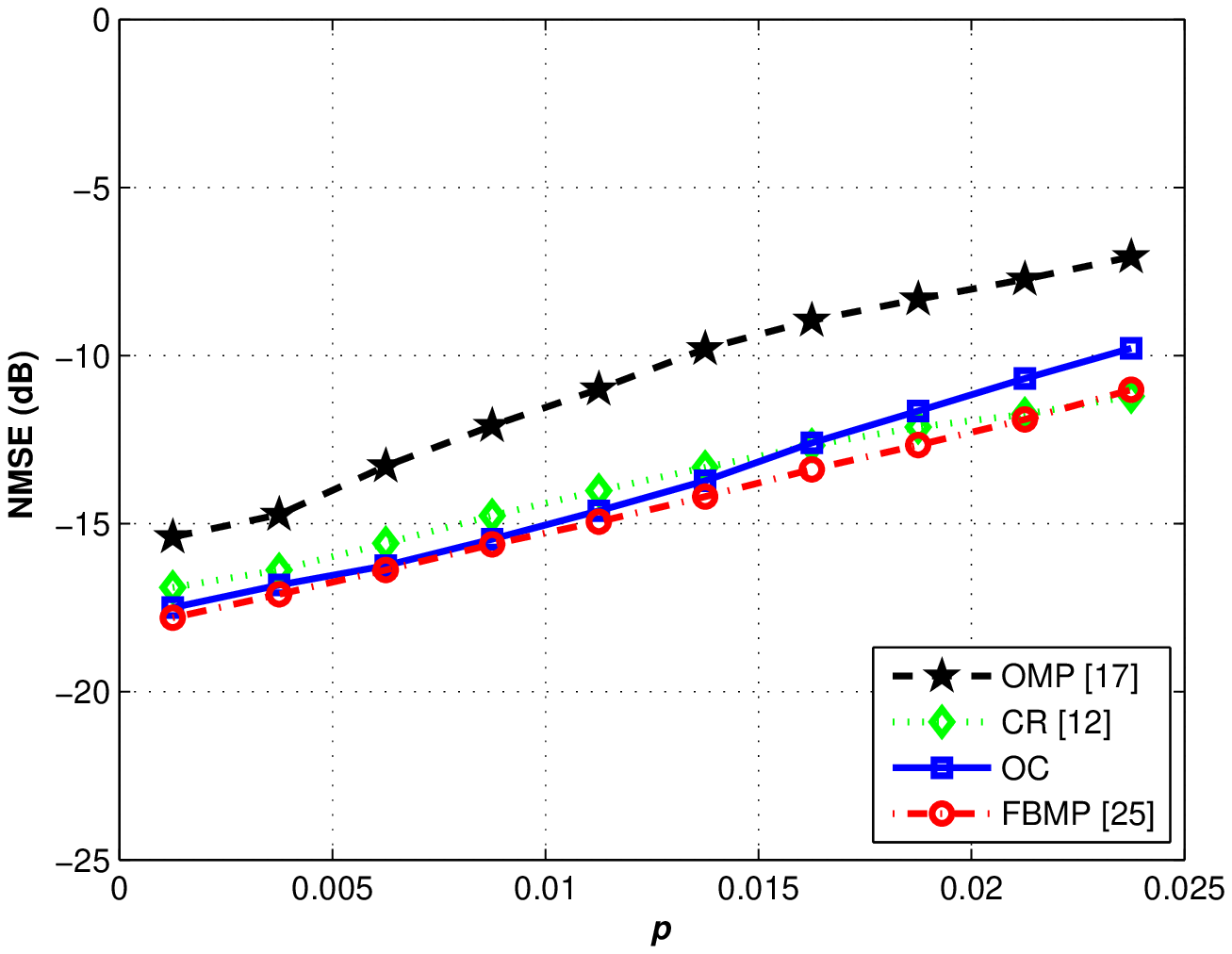}
\end{center}
\caption{\small NMSE vs $p$ for the DFT matrix and $\xv|\Sc$ Gaussian.} \label{nmse_p_DFT_Gaussian}
\end{minipage}
\begin{minipage}[htp]{.5\textwidth}
\begin{center}
\includegraphics[width=\textwidth]{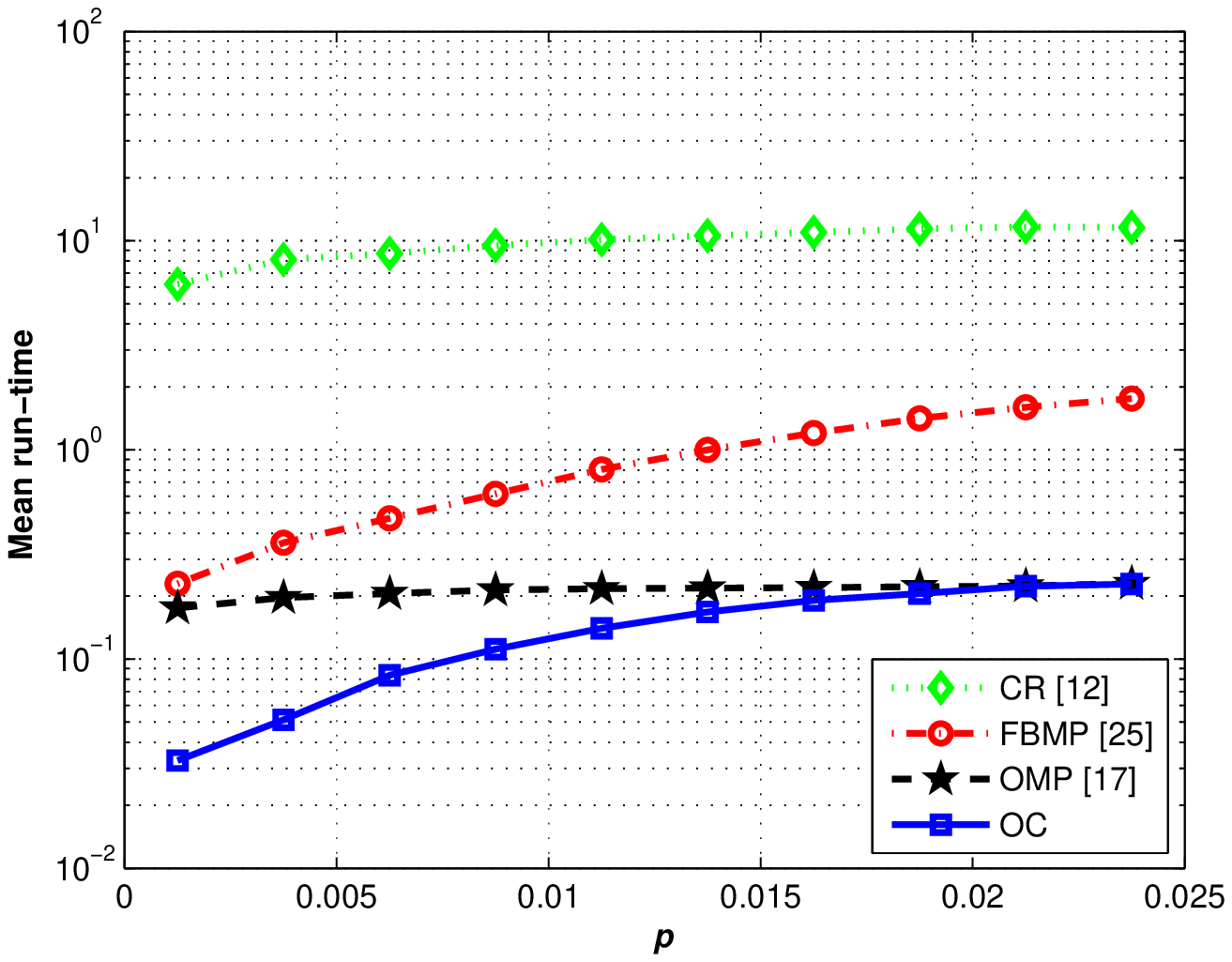}
\end{center}
\caption{\small Mean run-time for the DFT matrix and $\xv|\Sc$ Gaussian.} \label{runtime_p_DFT_Gaussian}
\end{minipage}
\begin{minipage}[htp]{.5\textwidth}
\begin{center}
\includegraphics[width=\textwidth]{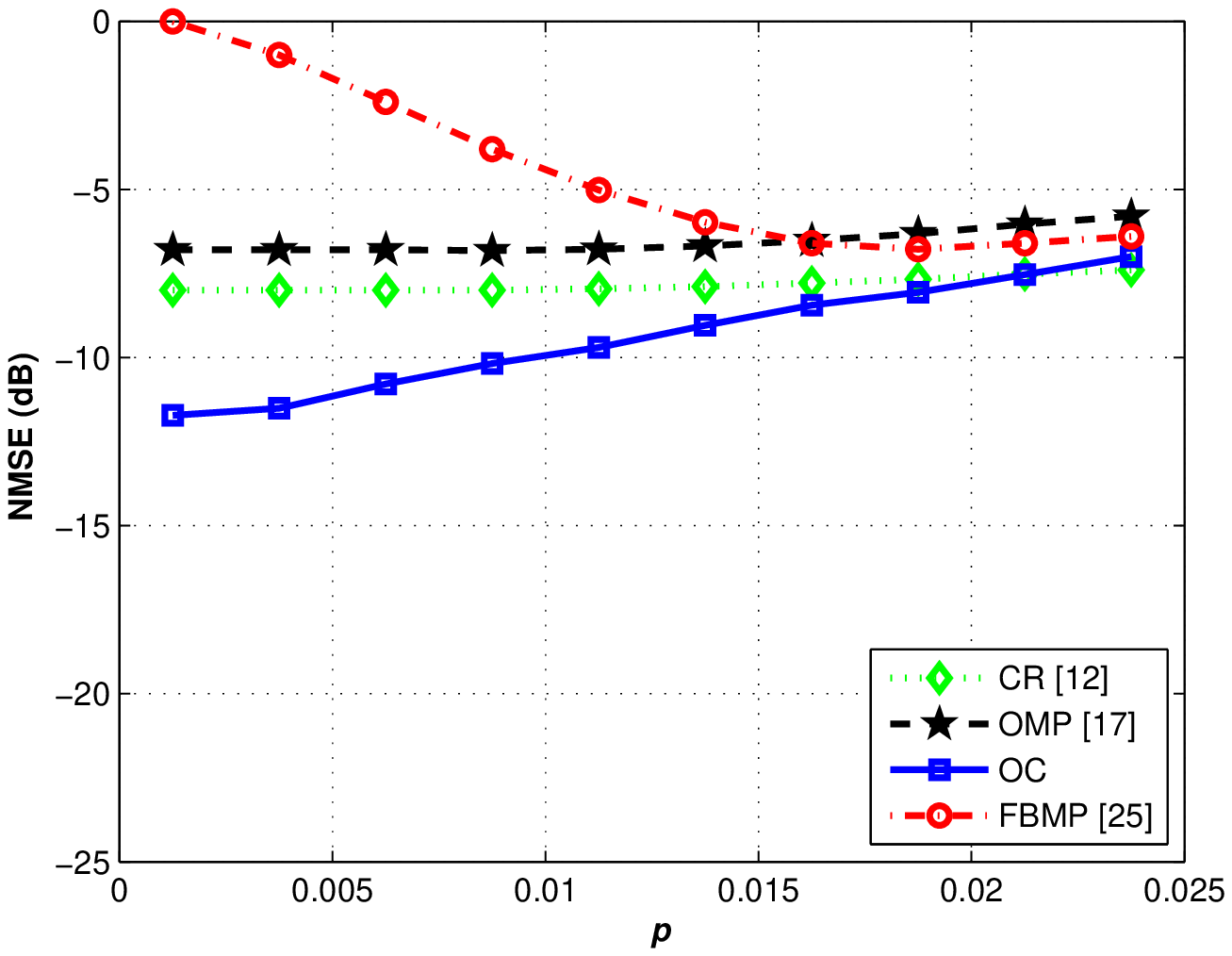}
\end{center}
\caption{\small NMSE vs $p$ for the DFT matrix and $\xv|\Sc$ unknown.} \label{nmse_p_DFT_nonGaussian}
\end{minipage}
\begin{minipage}[htp]{.5\textwidth}
\begin{center}
\includegraphics[width=\textwidth]{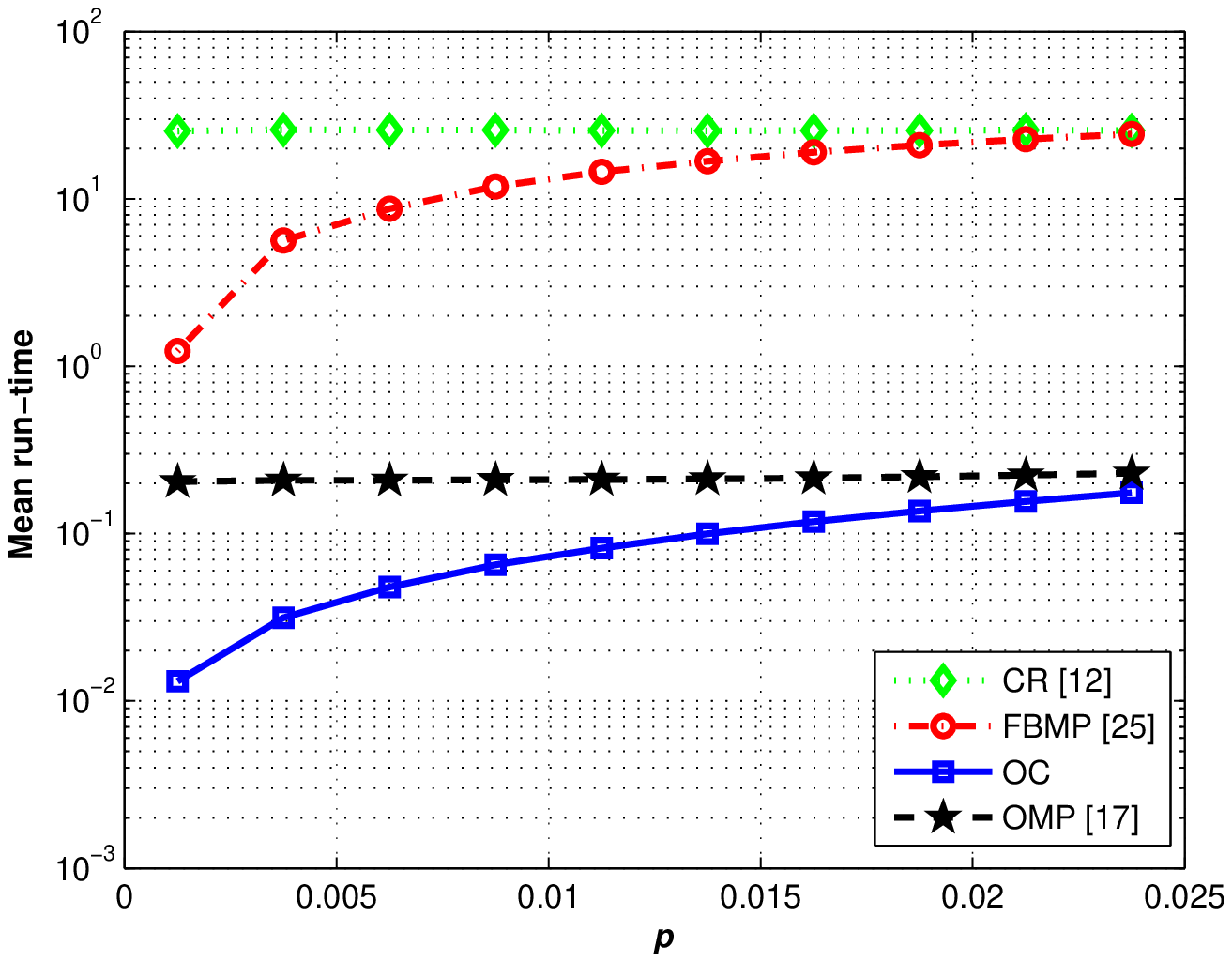}
\end{center}
\caption{\small Mean run-time for the DFT matrix and $\xv|\Sc$ unknown.} \label{runtime_p_DFT_nonGaussian}
\end{minipage}
\begin{minipage}[htp]{.5\textwidth}
\begin{center}
\includegraphics[width=\textwidth]{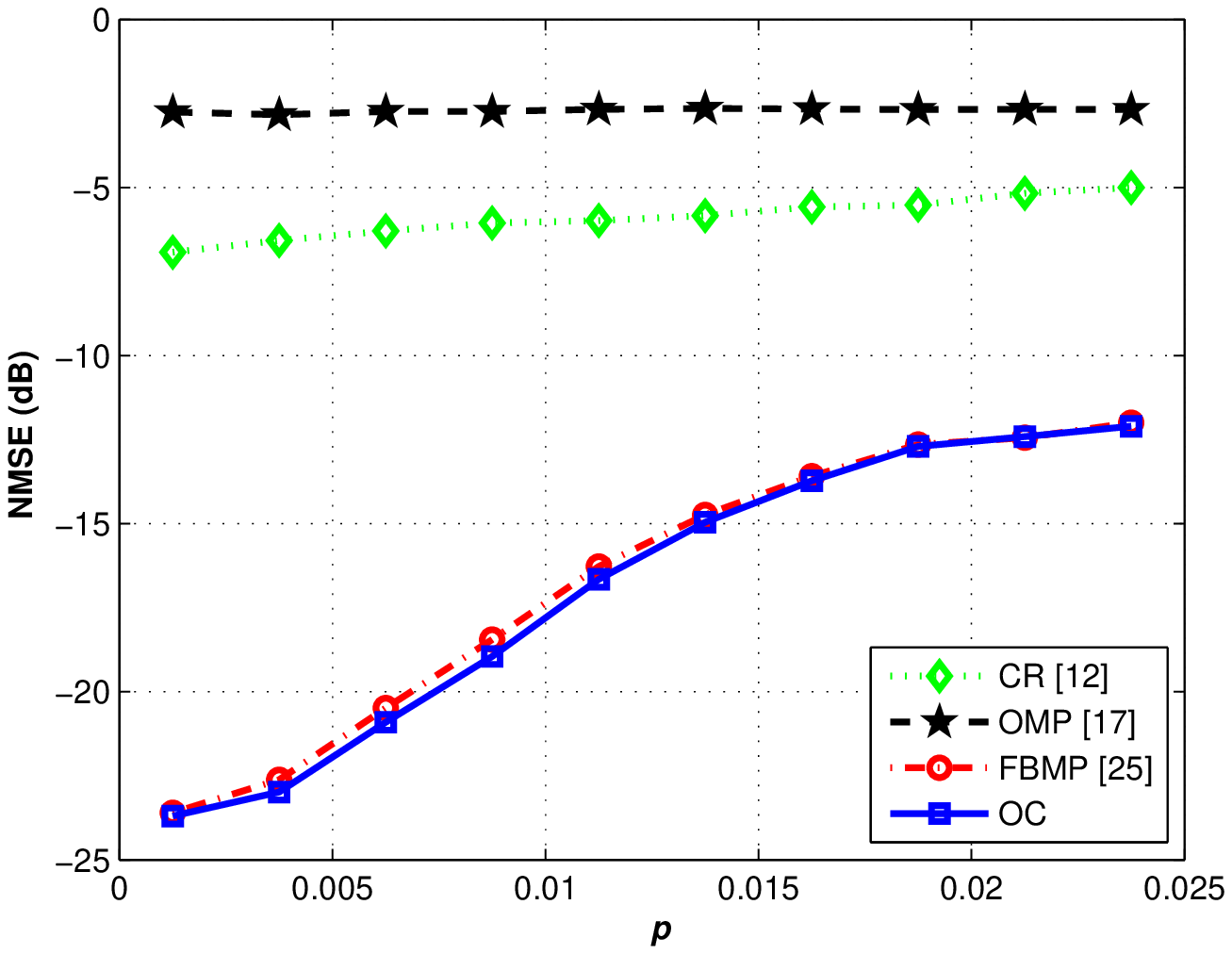}
\end{center}
\caption{\small NMSE vs $p$ for the Toeplitz matrix and $\xv|\Sc$ Gaussian.} \label{nmse_p_Toeplitz_Gaussian}
\end{minipage}
\begin{minipage}[htp]{.5\textwidth}
\begin{center}
\includegraphics[width=\textwidth]{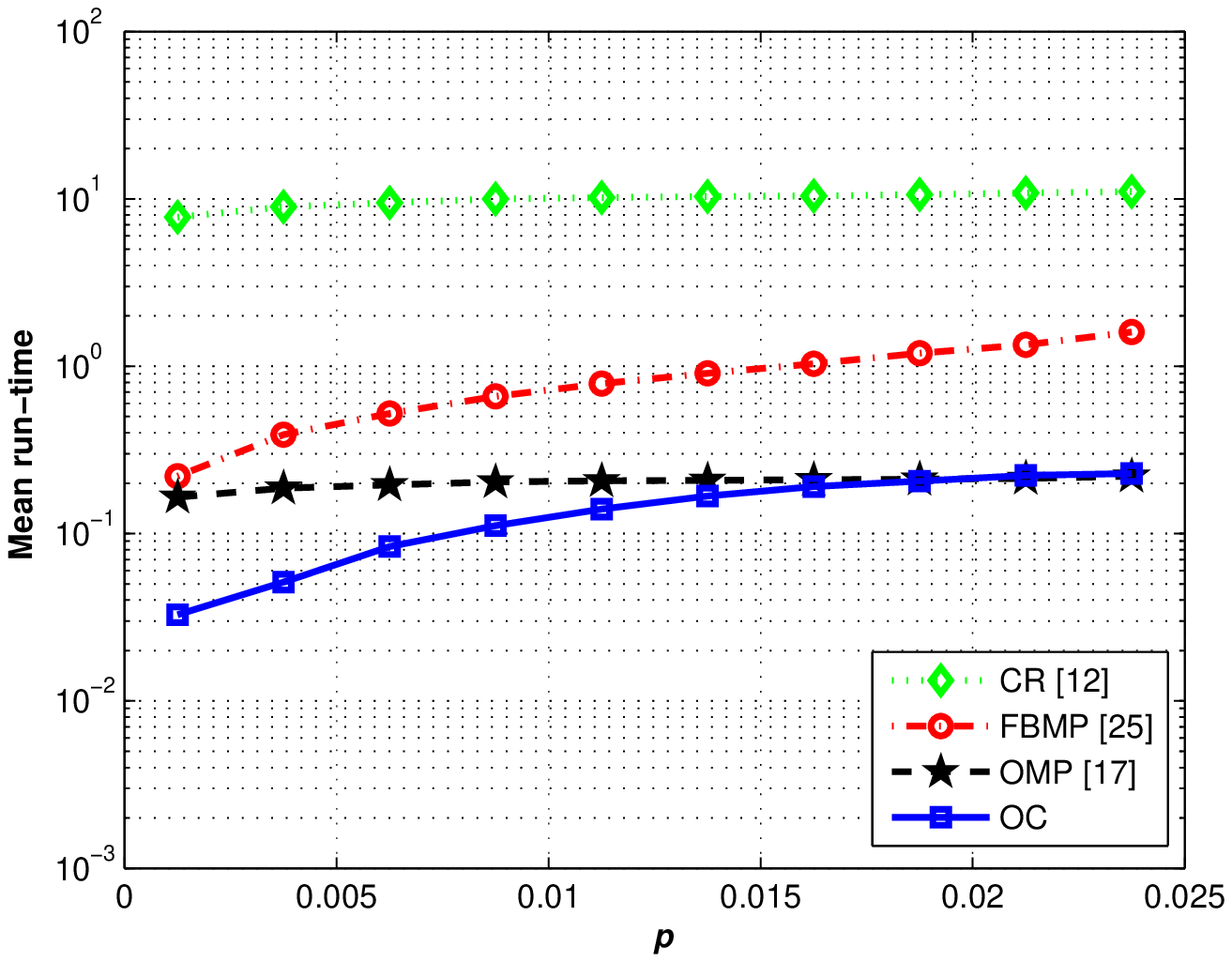}
\end{center}
\caption{\small Mean run-time for the Toeplitz matrix and $\xv|\Sc$ Gaussian.} \label{runtime_p_Toeplitz_Gaussian}
\end{minipage}
\end{figure}

\end{document}